\documentclass[12pt]{amsart}
\usepackage{amsfonts,latexsym,amsthm,amssymb,graphicx}
\usepackage[all]{xy}

\newtheorem{theorem}{Theorem}[section]
\newtheorem{lemma}[theorem]{Lemma}
\newtheorem{corol}[theorem]{Corollary}
\newtheorem{prop}[theorem]{Proposition}
\newtheorem{claim}[theorem]{Claim}

{\theoremstyle{definition} \newtheorem{defin}[theorem]{Definition}}
{\theoremstyle{remark} \newtheorem{remark}[theorem]{Remark}
\newtheorem{example}[theorem]{Example}}

\newcommand{\Cbb}{{\mathbb{C}}}

\newcommand{\Lbb}{{\mathbb{L}}}
\newcommand{\Pbb}{{\mathbb{P}}}
\newcommand{\Qbb}{{\mathbb{Q}}}
\newcommand{\Zbb}{{\mathbb{Z}}}

\newcommand{\Til}[1]{{\widetilde{#1}}}

\DeclareMathOperator{\codim}{codim}

\DeclareMathOperator{\Can}{Can}

\newcommand{\one}{1\hskip-3.5pt1}
\newcommand{\csm}{c_{\text{\rm SM}}}
\newcommand{\integ}[3][X]{\int_{#3}\one(#2)\,d\mathfrak c_#1}

\title{Modification systems and integration in their Chow groups}
\author{Paolo Aluffi}
\address{Max-Planck-Institut f\"ur Mathematik, Bonn, Deutschland}
\address{Math Dept, Florida State University, Tallahassee, Florida, U.S.A.}
\email{aluffi@math.fsu.edu}

\begin{document}

\begin{abstract}
We introduce a notion of integration on the category of proper
birational maps to a given variety $X$, with value in an associated
Chow group. Applications include new birational invariants; comparison
results for Chern classes and numbers of nonsingular birational
varieties; `stringy' Chern classes of singular varieties; and
a zeta function specializing to the topological zeta function.

In its simplest manifestation, the integral gives a new expression for
Chern-Schwartz-MacPherson classes of possibly singular varieties,
placing them into a context in which a `change-of-variable' formula
holds.
\end{abstract}

\maketitle

\setcounter{tocdepth}{1}
\tableofcontents

%%%

\newcommand{\cC}{{\mathcal C}}
\newcommand{\cD}{{\mathcal D}}
\newcommand{\cS}{{\mathcal S}}
\newcommand{\cZ}{{\mathcal Z}}
\newcommand{\cI}{{\mathcal I}}

\section{Introduction}\label{intro}

Our goal in this paper is the introduction of a technique for the
study of intersection theoretic invariants in the birational class of
a given variety $X$. Assuming resolution of singularities we introduce
an `integral', defined over the category $\cC_X$ of {\em
modifications\/} of (that is, proper and birational maps onto) $X$,
with values in the limit $A_*\cC_X$ of the Chow groups (with rational
coefficients) of the modifications. We will refer to the category
$\cC_X$ as the {\em modification system\/} of $X$. The push-forward of
the integral to a variety $Z$ is called its {\em manifestation\/} in
$Z$.  If two varieties $X$ and $Y$ admit a modification from a common
source then the corresponding systems $\cC_X$ and $\cC_Y$ can be
identified, in the sense that they share their main structures; for
example, their Chow groups are isomorphic. It is thus natural to look
for invariants of $X$ which can be expressed in terms of $\cC_X$: such
invariants have, via the identification of the latter with $\cC_Y$, a
natural manifestation as invariants of $Y$. Such invariants must really
be reflections of invariants of the birational class of $X$.

We find that {\em Chern classes\/} are such an invariant.  More
precisely, we can evaluate the total (homology) Chern class of the
tangent bundle of a nonsingular variety $X$ as the manifestation 
in $X$ of an integral $\integ{0}{\cC_X}$ over $\cC_X$
(Proposition~\ref{chernns}). By the mechanism described above, the
Chern class of $X$ has a natural manifestation in the Chow group of
any variety $Y$ sharing its modification system, even if no morphism
exists between $X$ and $Y$. 
(This operation does {\em not\/}
agree with the naive process of pulling back the class from $X$ to a
common modification and pushing-forward to $Y$; it is considerably
subtler.)

In the proper generality, the integral is defined over any
constructible subset $\cS$ of a system, and applies to a divisor $\cD$
of the system (these notions are defined in \S\S\ref{divis}
and~\ref{conssub}; constructible subsets and Cartier divisors of $X$
provide important examples). The `integral' is defined by
intersection-theoretic means in \S\ref{integral}, and is shown
(\S\ref{properties}) to satisfy additivity with respect to $\cS$,
a change-of-variables formula with respect to proper birational maps
$Y\to X$, and the `normalization' property with respect to Chern
classes mentioned above.

The mere existence of such an operation triggers several applications:
precise comparison results for Chern classes and numbers of birational
varieties, new birational invariants for nonsingular varieties,
`stringy' Chern classes of singular varieties, invariants of
singularities as contributions to a zeta function. These topics are
discussed in \S\ref{applications}.

For example, the relation with Chern classes implies immediately
that Chern numbers $c_1^i\cdot c_{n-i}$ of complete nonsingular
$n$-dimensional varieties in the same $K$-equivalence class
coincide. This is pointed out in \S\ref{otherK}; Theorem~\ref{spell}
and Corollary~\ref{elgen} give a blueprint to obtain identities
involving Chern numbers of birational varieties regardless of their
$K$-equivalence.

Stringy Chern classes may be defined by taking the lead from the
normalization mentioned above: since $\integ{0}{\cC_X}$ recovers
$c(TX)\cap [X]$ when $X$ is nonsingular, we can define the {\em
stringy\/} Chern class of $X$, for (possibly) singular varieties $X$,
to be the identity manifestation of $\integ{0}{\cC_X}$ (and again, the
definition of these objects through an integral provides us
simultaneously with a consistent choice of manifestations for all
varieties birational to $X$).

There is an important subtlety hidden in this definition: if $X$ is
singular, our integral depends on the consistent choice of a {\em
relative canonical divisor\/} for modifications $\pi: V\to X$ from a
nonsingular source, and there is more than one possible such
choice. Using K\"ahler differentials leads to what we call the
`$\Omega$~flavor' of the integral; this provides a stringy notion
which is defined for arbitrary singular varieties. In birational
geometry it is customary to use a notion arising from the double dual
of $\Omega^{\dim X}_X$, which leads to the `$\omega$~flavor' of our
integral. The corresponding stringy Chern class is defined for
$\Qbb$-Gorenstein varieties admitting at least one log resolution with
log discrepancies greater than zero, that is, with log terminal
singularities. If $X$ is complete, the degree of its $\omega$-stringy
Chern class agrees with Victor Batyrev's stringy Euler number
(\cite{MR2001j:14018}, with $\Delta_X=0$; the stringy Euler number for
a Kawamata pair $(X,\Delta_X)$ is the degree of
$\integ{-\Delta_X}{\cC_X}$). 

An example comparing the two flavors is given in \S\ref{compare}. The
delicate issue of the definition of the integral (and hence of the
stringy classes) over varieties with non log-terminal singularities is
discussed in \S\ref{negative}, but we fall short of a definition in
this case.

There is a different way to extract invariants for possibly singular
varieties $S$ from the integral defined here: one may embed $S$ in a
nonsingular ambient variety $X$, then take the integral of the divisor
$0$ {\em over the constructible set~$\cS$ of $\cC_X$ determined 
by~$S$.\/} We prove (Theorem~\ref{CSMconn}) that (in characteristic~0)
the identity manifestation of this integral agrees with a known
invariant, that is, the Chern-Schwartz-MacPherson class of $S$. In
fact, there is a very tight connection between the identity
manifestation of any integral and MacPherson's natural transformation,
discussed in \S\ref{CSM}. In particular, the stringy Chern classes of
a singular variety $X$ can be written as explicit linear combinations
of Chern-Schwartz-MacPherson classes of subvarieties of~$X$, as the
image through MacPherson's natural transformation of a specific
constructible function $I_X(0,\cC_X)$ on $X$. This constructible 
function appears to us as a much more basic invariant of $X$ than the
stringy Euler number or Chern class obtained from it; it would be
worth studying it further.

It should be noted that the theory of Chern-Schwartz-MacPherson
classes is {\em not\/} an ingredient of the integral introduced here;
thus, the result mentioned above provides us with a candidate for a
possible alternative treatment of Chern-Schwartz-MacPherson classes
(in the Chow group with rational coefficient), relying solely on
canonical resolution of singularities. The original definition of
MacPherson in \cite{MR50:13587} relied on transcendental invariants;
other approaches, such as the one in \cite{MR91h:14010}, appear to
require generic smoothness or be otherwise bound to characteristic
zero for more fundamental reasons. 

(Canonical) resolution of singularities is essential for our approach,
in two respects. First, the integral is defined by pushing forward a
weighted $\Qbb$-linear combination of classes, defined in a
nonsingular variety suitably `resolving' the data (\S\ref{intdef},
\S\ref{arbdef}); this exists by virtue of resolution of
singularities. Second, the key independence on the choice of resolving
variety (Claim~\ref{indepclaim}, Theorem~\ref{key}) is a technical
exercise relying on the factorization theorem of \cite{MR2003c:14016},
which uses resolution of singularities (cf.~Remark~(3), p.~533, in
\cite{MR2003c:14016}). Resolution of singularities appears to be the
only obstruction to extending the integral defined here to positive
characteristic.\vskip 6pt

Our guiding idea in this paper is an attempt to mimic some of the
formalism of {\em motivic integration\/} (as developed by Jan Denef
and Fran\c cois Loeser) in a setting carrying more naturally
intersection-theoretic information. While motivic integration is not
used in the paper, the reader may notice several points of contact; in
particular, the degree of our integral agrees with an expression which
arises naturally in that context (see especially
Remark~\ref{denloecon} and Claim~\ref{contact}). This observation was
at the root of our previous work along these lines,
\cite{math.AG/0401167}.

We found the excellent survey~\cite{math.AG/0401374} particularly
useful as a guide to motivic integration and its applications, and
refer to this reference in several place. We note that a different and
powerful approach to `motivic Chern classes' has appeared in the
recent work~\cite{math.AG/0405412}. Also, Shoji Yokura has examined
inverse systems of algebraic varieties for intersection-theoretic
purposes, in \cite{shoji}.

The most powerful approach to date to defining `stringy'
characteristic classes for singular variety appears to be the one
by Lev Borisov and Anatoly Libgober (\cite{MR1953295} and
\cite{math.AG/0206241}). While we have not checked this too carefully,
it is essentially inevitable that the ($\omega$ flavor of the) stringy
Chern class introduced here should agree with a suitable
specialization of the {\em elliptic orbifold class\/} introduced in
\cite{math.AG/0206241}, Definition~3.2. It is worth mentioning that
the factorization theorem of \cite{MR2003c:14016} is also at the root
of the approach of Borisov and Libgober, as it reduces the key
independence to the analysis of the behavior of the invariants through
blow-ups. Accordingly, our Theorem~\ref{key} could likely be derived
as a corollary of Theorem~3.5 in \cite{math.AG/0206241}. 

In this connection we would like to suggest that a notion of `elliptic
integration' could be defined using Borisov-Libgober's elliptic classes
analogously to what we have done here for Chern classes. It would be
very interesting if an analog of Theorem~\ref{CSMconn} were to hold
for such a notion, linking it to a functorial theory generalizing the
theory of Chern-Schwartz-MacPherson classes and MacPherson's natural
transformation (maybe the theory proposed in
\cite{math.AG/0405412}?). Lev Borisov informs me that he and Libgober
have entertained the idea of studying relations of their work with
other characteristic classes for singular varieties, such as
Chern-Mather classes.

The inverse limit of the modification system of a variety $X$ exists
as a ringed space, and Heisuke Hironaka calls it the {\em Zariski
space\/} of~$X$ (\cite{MR33:7333}, Chapter~0, \S6). We prefer to work
throughout with the inverse system of modifications, while
systematically taking limits of structures (Chow groups, etc.)
associated to it. As modification systems are closely related to
Hironaka's {\em vo\^ute \'etoil\'ee,\/} Matilde Marcolli has suggested
that we could name the operation defined here {\em celestial
integration.\/} We will show due restraint and use this poetic term
very sparingly. In any case we warmly thank her, as well as Ettore
Aldrovandi and J\"org Sch\"urmann, for extensive discussions on the
material presented here.

Thanks are also due to the Max-Planck-Institut f\"ur Mathematik in
Bonn for hospitality during the Summer of 2004, when a substantial part
of this work was done.

%%%

\section{Modification systems}\label{modsysts}

\subsection{}\label{prep}
We work over any algebraically closed field $k$ over which canonical
resolution of singularities \`a la Hironaka and the factorization
theorem of \cite{MR2003c:14016} hold.

In this section we introduce the objects (`modification systems') over
which the integration will take place. Roughly speaking, a
modification system consists of the collection of varieties mapping
properly and birationally onto a fixed variety $X$. A modification
system carries a number of natural structures, all obtained by taking
direct or inverse limits of corresponding structures on the individual
varieties in the collection: thus the {\em Chow group\/} of a
modification system is the inverse limit of the Chow groups (under
push-forward); a {\em divisor\/} is an element of the direct limit of
the groups of Cartier divisors (under pull-back); and so on.

The key advantage of working with modification system is that, thanks
to embedded resolution, the information carried by many such
structures can in fact be represented in terms of a divisor with
normal crossings and nonsingular components on one nonsingular
modification of $X$. Thus, for example, arbitrary subschemes of $X$
may be described in terms of divisors of the corresponding
modification system (cf.~Remark~\ref{remdivis}).

\subsection{}\label{modsystdef}
The main sense in which the quick description given in \S\ref{prep} is
imprecise is that working with the {\em varieties\/} mapping to $X$ is
inadequate; it is necessary to work with the proper, birational {\em
maps\/} themselves.

\begin{defin}
Let $X$ be an irreducible variety over $k$. 
The {\em modification system\/} $\cC_X$ of $X$ is the category whose
objects are the proper birational maps $\pi: V\to X$, and morphisms
$\alpha: \pi_1 \to \pi_2$ are commutative diagrams of proper
birational maps
$$\xymatrix{
V_1 \ar[rr]^\alpha \ar[dr]_{\pi_1} & & V_2 \ar[dl]^{\pi_2} \\
& X
}$$
\end{defin}

We often (but not always) denote by $V_\pi$ the source of $\pi$. If
$\pi$ and $\pi'$ are objects of a modification system $\cC_X$, we say
that $\pi$ {\em dominates\/} $\pi'$ if there is a morphism $\pi
\to\pi'$ in $\cC_X$.

\begin{lemma}\label{domin}
Every pair $\pi_1$, $\pi_2$ in $\cC_X$ is dominated by some object
$\pi$ of $\cC_X$. In fact, $\pi$ may be chosen so that its source $V_\pi$
is nonsingular, and its exceptional divisor has normal crossings and 
nonsingular components in $V_\pi$.
\end{lemma}

\begin{proof}
In fact, modification systems have products: if $\pi_1$, $\pi_2$ are in 
$\cC_X$, the component $W$ of $V_{\pi_1} \times_X V_{\pi_2}$ dominating
$V_{\pi_1}$ and $V_{\pi_2}$ makes the diagram
$$\xymatrix{
& W=V_\pi \ar[dl]_{\alpha_1} \ar[dr]^{\alpha_2} \ar[dd]^\pi \\
V_1 \ar[dr]_{\pi_1} & & V_2 \ar[dl]^{\pi_2}\\
& X
}$$
commute, with $\pi$ proper and birational and satisfying the evident
universal property.

This implies the first assertion. The second assertion follows from
embedded resolution of singularities.
\end{proof}

\subsection{}\label{equivdef}
Modification systems are inverse systems under the ordering $\pi \ge \pi' \iff$ 
$\pi$ dominates $\pi'$. We want to identify systems `with the same limit'; in 
our context, this translates into the following.
\begin{defin}
Two systems $\cC_X$, $\cC_Y$ are {\em equivalent\/} if there are
objects $\pi_X$ in $\cC_X$ and $\pi_Y$ in $\cC_Y$ with isomorphic source;
that is, if there exists a variety $V$ and proper birational maps
$\pi_X:V\to X$ and $\pi_Y:V\to Y$: 
$$\xymatrix{
& V \ar[dl]_{\pi_X} \ar[dr]^{\pi_Y}\\
X && Y
}$$
\end{defin}

The transitivity of this notion follows from Lemma~\ref{domin}.

\begin{lemma}\label{immed}
If $\cC_X$ and $\cC_Y$ are equivalent, then every object $\rho$ of
$\cC_Y$ is dominated by an object $\sigma$ whose source is the source
of an object of $\cC_X$.
\end{lemma} 

\begin{proof}
This follows immediately from Lemma~\ref{domin}:
since $\cC_X$ and $\cC_Y$ are equivalent, there are objects $\pi_X$, 
$\pi_Y$ in $\cC_X$, $\cC_Y$ resp., with a common source $V$; any object
$\sigma$ dominating both $\pi_Y$ and $\rho$:
$$\xymatrix{
& & W \ar[dl] \ar[dr]  \ar@/_1pc/[ddrr]_{\sigma} \\
& V\ar[dl]^{\pi_X} \ar[drrr]_{\pi_Y} & & {V_\rho} \ar[dr]^\rho\\
X & & & & Y
}$$
satisfies the requirement.
\end{proof}

\subsection{}\label{copy}
A proper birational morphism $\rho: Y \to X$ determines a covariant
functor $\cC_Y \to \cC_X$ by composition: $\pi\mapsto \rho\circ\pi$;
$\cC_X$ and $\cC_Y$ are then trivially equivalent systems. 

Via this functor, the category $\cC_Y$ is realized as a {\em full\/}
subcategory of $\cC_X$. Indeed, every morphism $\rho\circ\pi_1 
\to\rho\circ\pi_2$ in $\cC_X$ is given by a proper birational map
$\alpha$ such that $\rho\circ\pi_1=\rho\circ\pi_2\circ\alpha$; but
then $\pi_1$ and $\pi_2\circ\alpha$ agree on a nonempty open 
set (on which they are isomorphisms), hence they agree everywhere. 
In other words, every morphism $\rho\circ\pi_1 \to \rho\circ\pi_2$
is induced by a morphism $\pi_1 \to \pi_2$ in $\cC_Y$.

In particular, any object $\rho: V\to X$ of $\cC_X$ determines a
modification system equivalent to $\cC_X$ and equal to a `copy' of
$\cC_V$; we will denote this system by $\cC_\rho$. Two systems
$\cC_X$, $\cC_Y$ are equivalent if and only if they both contain a
copy of a third system~$\cC_V$.

As mentioned in \S\ref{prep}, standard structures can be defined on
modification systems by taking limits of the same structures on the
(sources of the) objects in the system. The foregoing considerations
imply that any functorial structure defined in this fashion will be
preserved under equivalence of systems.

We will explicitly need very few such structures: Chow groups,
divisors, constructible sets. These are presented in the next few
sections.

\subsection{}\label{Chow}
For $V$ a variety, we denote by $A_*V$ the Chow group of $V$,
{\em tensored with $\Qbb$\/} (rational coefficients appear to be 
necessary for the integral introduced in \S\ref{integral}). The 
modification system $\cC_X$ determines the inverse system of Chow
groups $A_*V_\pi$ (under push-forwards), as $\pi$ ranges over the 
objects of $\cC_X$.

\begin{defin}
The {\em Chow group\/} of $\cC_X$ is the limit
$$A_*\cC_X:=\varprojlim_{\pi\in Ob(\cC_X)} A_*V_\pi\quad.$$
\end{defin}

That is, an element of the Chow group of $\cC_X$ is the choice of
a class $a_\pi$ in $A_*V_\pi$ for each object $\pi$ of $\cC_X$, with the 
condition that if $\alpha:\pi_1 \to \pi_2$ is a morphism then 
$\alpha_*(a_{\pi_1})=a_{\pi_2}$. We say that $a_\pi$ is the {\em 
$\pi$-manifestation\/} of this element of $A_*\cC_X$.

If $X$ is complete, we may define the {\em degree\/} of a class $a\in
A_*\cC_X$,
$$\deg a$$
as the degree of (the zero-dimensional component of) any manifestation
$A_\pi\in A_*V_\pi$ of $a$; indeed, this is clearly independent of~$\pi$.

Next, we observe that equivalent systems have isomorphic Chow groups.

\begin{lemma}\label{isoeq}
For all objects $\rho$ of $\cC_X$, $A_*\cC_\rho\cong A_*\cC_X$.
\end{lemma}

\begin{proof}
If $\sigma$ is an object of $\cC_\rho$ then $\rho\circ\sigma$ is an
object of $\cC_X$. The identity homomorphism $A_*V_{\rho\circ\sigma} \to
A_*V_\sigma$ induces compatible homomorphisms $A_*\cC_X\to
A_*V_\sigma$ for all $\sigma$; and hence a homomorphism $A_*\cC_X\to
A_*\cC_\rho$. 

Lemma~\ref{domin} gives the inverse homomorphism. More
precisely, for $\pi$ in $\cC_X$ there exists a $\pi'$ dominating both
$\pi$ and $\rho$:
$$\xymatrix{
& W \ar[dl] \ar[dd]^{\pi'} \ar[dr]\\
V_\rho \ar[dr]_\rho & & V_\pi \ar[dl]^\pi\\
& X
}$$
Push-forwards yield a homomorphism $A_*\cC_\rho \to A_*W \to
A_*V_\pi$, which is easily seen to be independent of the chosen $\pi'$
dominating $\pi$ and $\rho$, and to satisfy the necessary
compatibility, giving a homomorphism $A_*\cC_\rho \to A_*\cC_X$ by 
the universal property of the inverse limit.

Checking that the two compositions are the identity is equally
straightforward.
\end{proof}

\begin{corol}\label{isoeq2}
Equivalent systems have isomorphic Chow groups.
\end{corol}

\begin{proof}
Indeed, equivalent systems $\cC_X$, $\cC_Y$ contain a copy of a third
system $\cC_V$, cf.~\S\ref{copy}. By Lemma~\ref{isoeq}, $A_*\cC_X\cong
A_*\cC_V\cong A_*\cC_Y$.
\end{proof}

Note however that these isomorphisms are not canonical, as they 
depend on the choice of a common source of objects in $\cC_X$ and
$\cC_Y$. A canonical isomorphism is available if the common source is
chosen; for example, for a given proper birational map $\pi: Y\to X$.

\subsection{}\label{divis}
Divisors are arguably the most important structure within a
modification system. The group of divisors of $\cC_X$ is the direct
limit of the groups of Cartier divisor of sources of objects of $\cC_X$.
Explicitly:

\begin{defin}\label{divisdef}
Let $\cC_X$ be a modification system. A {\em divisor\/} $\cD$ of 
$\cC_X$ is a pair $(\pi,D_\pi)$, where $\pi$ is an object of $\cC_X$
and $D_\pi$ is a Cartier divisor on the source $V_\pi$ of $\pi$,
modulo the equivalence relation:
\begin{multline*}
(\pi',D_{\pi'}) \sim (\pi'',D_{\pi"}) \iff \text{for all 
$\pi=\pi'\circ\alpha'=\pi''\circ\alpha''$ in $\cC_X$}\\
\text{dominating $\pi'$ and $\pi''$, ${\alpha'}^*D_{\pi'}
={\alpha''}^*D_{\pi"}$ in $V_\pi$.}
\end{multline*}
\end{defin}
That is: $D_{\pi'}$ and $D_{\pi"}$ determine the same divisor of 
$\cC_X$ if they agree after pull-backs.

\begin{remark}\label{remdivis}
\begin{itemize}

\item Divisors of $X$ determine divisors of $\cC_X$ by pull-backs. In
fact, every proper subscheme $Z$ of $X$ determines a divisor of
$\cC_X$, namely $(\pi,E)$, where $\pi:\Til X \to X$ is the
blow-up along $Z$, and $E$ is the exceptional divisor.

\item In particular, if $X$ is smooth then every canonical divisor
of $X$ determines a divisor of $\cC_X$. If $\cC_X$ and $\cC_Y$ are
equivalent, the corresponding canonical divisors need not agree; if
they do, then $X$ and $Y$ are in the same {\em $K$-equivalence
class.\/} This notion has been thoroughly studied by Chin-Lung Wang
(for $\Qbb$-Gorenstein varieties), see \cite{MR2003j:14015}. 

\item As we are assuming that resolution of singularities hold, 
every divisor of $\cC_X$ admits a representation $(\pi,D)$
in which the source $V_\pi$ of $\pi$ is nonsingular, and~$D$ is
supported on a divisor with normal crossings and nonsingular 
components.

\item Divisors act on the Chow group: if $\cD$ is a divisor and $a\in
  A_*\cC_X$, for every object $\rho$ we can find a dominating
  $\pi=\rho\circ\alpha$ such that $\cD$ is represented by
  $(\pi,D_\pi)$, and set
$$(\cD\cdot a)_\rho:= \alpha_*(D_\pi\cdot a_\pi)\in A_*V_\rho\quad;$$
the projection formula guarantees that this is independent of the
chosen $\pi$, and that the resulting classes satisfy the compatibility
needed to define an element $\cD\cdot a\in A_*\cC_X$.

\item Equivalent systems have isomorphic divisor groups. Indeed, 
if $\cC_X$ and $\cC_Y$ are equivalent, and $(\pi, D_\pi)$ represents 
a divisor on $\cC_X$, then (by Lemma~\ref{immed}) we may assume
that $V_\pi$ is a source of an object of $\cC_Y$, so that $D_\pi$
represents a divisor on $\cC_Y$; and, conversely, divisors of $\cC_Y$
determine divisors of $\cC_X$.

As in \S\ref{Chow}, these identifications are canonical once a common
source is chosen.
\end{itemize}
\end{remark}

\subsection{}\label{conssub}
A second main character in the definition of the integral is the notion
of constructible subset.
\begin{defin}
A {\em constructible subset\/} $\cS$ of a modification system $\cC_X$ is
a pair $(\pi,S_\pi)$, where $\pi$ is an object of $\cC_X$ and $S_\pi$ is a 
constructible subset of $V_\pi$, modulo the equivalence relation:
\begin{multline*}
(\pi',S_{\pi'}) \cong (\pi'',S_{\pi''}) \iff \text{for all
$\pi=\pi'\circ\alpha' =\pi''\circ\alpha''$ in $\cC_X$}\\
\text{dominating $\pi'$ and $\pi''$,
${\alpha'}^{-1}(S_{\pi'})={\alpha''}^{-1}(S_{\pi''})$ in $V_\pi$.}
\end{multline*}
\end{defin}

That is: $S'$ and $S''$ determine the same subset of $\cC_X$ if they
agree after preimages through the system.

\begin{remark}
\begin{itemize}\label{remconssub}
\item As for divisors, note that the constructible subsets of a system
  are in one-to-one correspondence with constructible subsets of
  equivalent systems.
\item By resolution of singularities, every constructible subset of a
  system may be represented by a pair $(\pi,S_\pi)$ where $\pi$ has
  nonsingular source, and $S_\pi$ is obtained by taking unions and
  complements of nonsingular hypersurfaces meeting with normal
  crossings.
\item One can take unions or intersections of constructible subsets,
 by performing these operations in a $V_\pi$ in which all terms admit
 representatives.
\end{itemize}
\end{remark}

We say that $\cS$ is {\em closed, locally closed, etc.\/} if it is
represented by a closed, locally closed, etc.~set.

\subsection{}
More structures could be considered easily, such as maps from a
modification system to a variety, or constructible functions on a
modification system, etc. We will only have fleeting encounters with
such notions, and the reader should have no difficulties filling in
appropriate definitions as needed.

%%%

\newcommand{\Jbb}{{\mathbb J}}
\newcommand{\cT}{{\mathcal T}}

\section{Definition of the integral}\label{integral}

\subsection{}
{\em Caveat on terminology.\/} In this section we introduce the
`integration' operation on a modification system. As the reader may
now expect the appearance of a measure, a special class of
functions, and other ingredients of a good theory of integration,
we hasten to warn that none will be given here.

Our guiding idea comes from {\em motivic integration\/} (see e.g.,
Eduard Looijenga's Bourbaki survey, \cite{MR2003k:14010}). A motivic
integral is obtained from a suitable measure on the arc space of a
variety $X$; the objects to which it is applied are in the form
$\Lbb^{-\alpha}$, where $\alpha$ is a constructible function on the
space. In practice, most applications are to the case in which
$\alpha$ is the order function of a divisor on $X$; and resolution of
singularities allows one to further restrict attention to functions
arising from a divisor with normal crossings and nonsingular
components.  In this case the integral can be computed explicitly,
bypassing motivic measure entirely: see for example Alastair Craw's
``user-friendly formula'' (Theorem 2.15 in \cite{math.AG/9911179} or
\S3.6 in \cite{math.AG/0401374}).

Our definition is motivated by such formulas. As pointed out in
\S\ref{modsysts}, most structures of importance in the context of
modification systems are encoded by normal crossing divisors,
so we feel free to cut the middle man and offer the `user-friendly'
analog as our definition. It would be interesting to interpret our
definition in terms of a measure (be it on the modification system, or
perhaps on its arc space), but this does not appear to be necessary
for applications.

One advantage of this approach is that it comes with a built-in
change-of-variables formula. We pay the price for this benefit by
having to prove independence of the choice of representative for the
divisor. In motivic integration the analogous formula is proved to
agree with the definition based on a measure, so the corresponding
independence is automatic; while the change-of-variables formula
requires an argument.

\subsection{}\label{relcandiv}
Let $X$ be an irreducible variety, let $\cD$ be a divisor in the
modification system $\cC_X$ of $X$ (cf.~\S\ref{divis}), and let $\cS$
be a constructible subset of $\cC_X$ (cf.~\S\ref{conssub}). The rest
of this section is devoted to the definition of an element
$$\integ{\cD}{\cS} \in A_*\cC_X\quad.$$
While the definition is rather transparent, proving that it does
not depend on the various choices will require a certain amount of
technical work. Simple properties of this definition, and
applications, will be discussed in later sections and will not involve
the more technical material in the present one.

We have to define an element
$$\integ{\cD}{\cS}$$
of $A_*\cC_X$, and this is equivalent to defining all of its
manifestations in $A_*V_\pi$, as $\pi$ ranges over the objects of
$\cC_X$.  In order to streamline the exposition, we will begin by
assuming that $\cS$ is closed, and that $\cD$, $\cS$, and the {\em
relative canonical divisor\/} are in a particularly favorable position
in $V_\pi$; the definition in this case is given in
\S\ref{intdef}. The definition for all objects of $\cC_X$ is given in
\S\ref{arbdef}, and the extension to constructible subsets $\cS$ is
completed in~\S\ref{constr}.

The notion of relative canonical divisor requires a discussion. If
$\pi: V \to X$ is a birational morphism of nonsingular varieties, we
will denote by $K_\pi$ the divisor of the jacobian of~$\pi$; so if
$K_X$ is a canonical divisor of $X$, then $\pi^*K_X+K_\pi$ is a
canonical divisor of $V$. This notion behaves well with respect to
composition of maps, in the sense that if $\alpha: W \to V$ is a
birational morphism of nonsingular varieties, then
\begin{equation*}\tag{*}
K_{\pi\circ\alpha}=K_\alpha+\pi^* K_\pi\quad.
\end{equation*}

We must call the attention of the reader to the fact that there are
several ways to extend this notion to the case in which $X$ may be
{\em singular.\/} Necessary requirements from our viewpoint are that
\begin{enumerate}
\item\label{agreens} the notion agrees with the one recalled above in
the nonsingular case;
\item\label{getdivnc} if $\pi: V \to X$ is a proper birational
morphism, with $V$ nonsingular, then there exists a $\pi'$ dominating
$\pi$ and such that $K_{\pi'}$ is a divisor with normal crossings and
nonsingular components in~$V_{\pi'}$; and 
\item\label{funrcd} (*) holds for 
$\xymatrix@1{W \ar[r]^\alpha & V \ar[r]^\pi & X}$,
with $V$ and $W$ nonsingular, assuming $K_\pi$ is a divisor.
\end{enumerate}
Note that we are not requiring that $K_\pi$ be a divisor {\em for
all\/} $\pi$. This may lead to some confusion as we will nevertheless
stubbornly refer to any $K_\pi$ as a relative canonical {\em
divisor,\/} since it determines a divisor in the system
(cf.~Remark~\ref{remdivis}).

One simple possibility, which has the advantage of working without
further assumptions on $X$, is the following. If $\pi: V\to X$ is a
birational morphism of $n$-dimensional varieties, there is an induced
morphism of sheaves of K\"ahler differentials
$$\pi^*\Omega^n_X \to \Omega^n_V\quad.$$
If $V$ is nonsingular, so that $\Omega^n_V$ is locally free, the image
of this morphism may be written as $\Omega^n_V\otimes \mathcal I$ for
an $\mathcal O_V$-ideal sheaf $\mathcal I$. By composing with a
sequence of blow-ups $\rho:V' \to V$ we may ensure that the
ideal $\mathcal I'$ in $\mathcal O_{V'}$ corresponding to
$\pi'=\pi\circ\rho$ is principal, thus $\mathcal I'=\mathcal
O_{V'}(-K_{\pi'})$ for a Cartier divisor $K_{\pi'}$ on $V'$. In fact,
by the same token $K_{\pi'}$ may be assumed to be a divisor with
normal crossings and nonsingular component, as promised.

This notion of $K_\pi$ gives what we call the `$\Omega$ flavor' of the
integral. A second possibility, which is more natural in the context of
birational geometry, will be mentioned later (\S\ref{singclass}) and
gives the `$\omega$ flavor'. 

It should be noted that, for $X$ singular, the value of the integral
will depend on the chosen notion of relative canonical divisor,
cf.~Example~\ref{compare}; but the basic set-up and properties do not
depend on this choice, so we will not dwell further on this important
point until later sections. By requirement (\ref{agreens}) above,
the integral is univocally determined if the base~$X$ is nonsingular.

\subsection{}\label{intdef}
Let $\pi$ be an object of $\cC_X$, $K_\pi$ be the relative canonical
divisor of $\pi$, and let $\cD=(\pi,D)$, $\cS=(\pi,S)$. {\em We assume
that $V_\pi$ is nonsingular, $K_\pi$ is a divisor, and there exists a
divisor $E$ of $V_\pi$, with normal crossings and nonsingular
components $E_j$, $j\in J$, such that 
$$D+K_\pi = \sum_{j\in J} m_j E_j\quad,$$
with $m_j>-1$; further, we assume that $S$ is $V_\pi$ or the union of
a collection of components $\{E_\ell\}_{\ell\in L}$ of $E$, and we let
$\Jbb_S$ be the whole family of subsets of $J$, if $S=V_\pi$,
or the subfamily of those subsets which meet $L$ if $S=\cup_{\ell\in
L} E_\ell$.\/}

In other words, we essentially require $\pi$ to be a `log resolution'
for the relevant data. We will say that $\pi$ (or, more loosely,
$V_\pi$) {\em resolves\/} $\cD$, $\cS$ if it satisfies these assumptions. 

\begin{remark}
The condition that no $m_j$ be $\le -1$ is undesirable, but necessary
for the present set-up. Ideally one would like to replace this with
the weaker request that $m_j$ be $\ne -1$ for all $j$, which suffices
for the expression in Definition~\ref{mandef} to make sense. But this
would come at the price of having certain manifestations of the
integral remain undefined; more importantly, the argument given here
does not suffice to prove that this would lead to consistent
definitions in the presence of multiplicities~$\le -1$.

This issue is discussed further and illustrated with an example in
\S\ref{negative}.
\end{remark}

\begin{defin}\label{mandef}
If $\pi$ resolves $\cD$, $\cS$, then the manifestation of $\integ{\cD}{\cS}$ 
in $A_*V_\pi$ is defined to be
$$\left(\integ{\cD}{\cS}\right)_{\pi}:=
c(TV_\pi(-\log E))\cdot\left(\sum_{I\in \Jbb_S}\prod_{i\in I}\frac{E_i}
{1+m_i}\right)\cap [V_\pi]$$
\end{defin}
Here $TV_\pi(-\log E)$ is the dual of the bundle of differential
forms with logarithmic poles along $E$. Its Chern class serves as a
shorthand for a longer expression, cf.~Lemma~\ref{necfacts},
(\ref{log}).

\begin{remark}\label{support}
\begin{itemize}
\item The expression given in Definition~\ref{mandef} is supported on
$S$; it should be viewed as the push-forward to $V_\pi$ of an element
of $A_*S$.
\item If $E_j$ is a component which does not belong to $S$, and for
which $m_j=0$, then the given expression is independent of whether
$E_j$ is counted or not in~$E$ (exercise!).
\end{itemize}
\end{remark}

\subsection{}\label{altexp}
The expression given in Definition~\ref{mandef} can be written in 
several alternative ways, some of which are rather suggestive, and
sometimes easier to apply.

For example, for $I\subset J$ let $E_I$ equal the intersection $\cap_{i\in 
I} E_i$; this is a nonsingular subvariety of $V_\pi$ since $E$ is a 
divisor with normal crossings and nonsingular components. Also, denote 
by $E^{\overline I}$ be the divisor $\sum_{i\not\in I} E_i$, as well as
its restriction to subvarieties of $V_\pi$; note that $E^{\overline I}$ 
intersects the subvariety $E_I$ along a divisor with normal crossings and 
nonsingular components.

Then the given expression of $\integ{\cD}{\cS}$ in $A_*V_\pi$ equals
$$\sum_{I\in \Jbb_S} \frac{c(TE_I(-\log E^{\overline I}))\cap [E_I]}
{\prod_{i\in I} (1+m_i)}\quad.$$

If $S=V_\pi$, so that $\Jbb_S$ is the whole family of subsets of $J$, then 
trivial manipulations show that the class equals a linear combination of
the Chern classes of the subvarieties~$E_I$:
$$\sum_{I\subset J} (-1)^{|I|}\prod_{i\in I} \frac{m_i}{1+m_i} \cdot
c(TE_I)\cap [E_I]\quad.$$
Equally trivial manipulations show that the same class can be written
as a weighted average of $\log$-twists of the Chern class of $V_\pi$:
$$\frac 1{\prod_{j\in J} (1+m_j)} \sum_{I\subset J} m_I\cdot 
c(TV_\pi(-\log E^I))\cap [V_\pi]$$
where $m_I=\prod_{i\in I}m_i$, and $E^I=\sum_{i\in I} E_i$; as 
$\sum_{I\subset J} m_I=\prod_{j\in J} (1+m_j)$, the top-dimensional term 
in this expression equals $[V_\pi]$, as it should.

\subsection{}\label{arbdef}
We have to define the manifestation of $\integ{\cD}{\cS}$ in
$A_*V_\rho$, for arbitrary objects $\rho$ of $\cC_X$.  By
principalization and embedded resolution of singularities, for any
divisor $\cD$ and closed subset $\cS$, every object $\rho$ in $\cC_X$
is dominated by an object $\pi$ resolving $\cD$, $\cS$.

\begin{defin}\label{closedS}
Let $\cD$ be a divisor of $\cC_X$, and let $\cS$ be a closed subset
of $\cC_X$. For arbitrary $\rho$ in $\cC_X$, let $\alpha$ be a proper
birational map such that $\pi=\rho\circ\alpha$ resolves $\cD$, $\cS$. 
Then the manifestation of $\integ{\cD}{\cS}$ in $A_*V_\rho$ is 
defined to be
$$\left(\integ{\cD}{\cS}\right)_\rho:=
\alpha_*\left(\integ{\cD}{\cS}\right)_\pi\quad.$$
\end{defin}

Of course we have to prove that this expression does not depend on the
choice of a resolving $\pi$ dominating $\rho$, and we note that this
will also immediately imply that the manifestations do define an
element of the inverse limit $A_*\cC_X$.

To establish the independence on the choice, we have to prove
that if $\pi_1=\rho\circ\alpha_1$ and $\pi_2=\rho\circ\alpha_2$
both resolve $\cD$, $\cS$, then the two push-forwards
${\alpha_i}_*\left(\integ{\cD}{\cS}\right)_{\pi_i}$
to $A_*V_\rho$, $i=1,2$, coincide. This is our next task.

By the factorization theorem of \cite{MR2003c:14016}, there exists
an object $\pi$ dominating both $\pi_1$ and $\pi_2$:
$$\xymatrix{
& V_\pi \ar[dl]_{\beta_1} \ar[dr]^{\beta_2}\\
V_{\pi_1} \ar[dr]_{\alpha_1} & & V_{\pi_2} \ar[dl]^{\alpha_2}\\
& V_\rho \ar[d]_{\rho} \\
& X
}$$
and such that $\beta_i$ decomposes as a sequence of maps
$$\xymatrix{
V_\pi \ar@/_1pc/[rr]_{\beta_i} \ar[r]^{\gamma_1} & \cdots 
\ar[r]^{\gamma_r} & V_{\pi_i}
}$$ 
with each $\gamma_k$ a sequence of blow-ups followed by a sequence
of blow-downs. Further, the centers of these blow-ups may be chosen
to intersect the relevant divisors with normal crossings.

We should note that \cite{MR2003c:14016} assumes the varieties
to be complete; we can reduce to this case by working in the modification
system $\cC_{\overline X}$ of a completion $\overline X$ of $X$, and then
taking the inverse image of $X$ throughout the system.

\begin{claim}\label{indepclaim}
If $\pi_1$ and $\pi_2$ resolve $\cD$, $\cS$, then so do all the 
intermediate stages in the resolution. In particular, so does $\pi$;
further,
$$\left(\integ{\cD}{\cS}\right)_{\pi_i}
={\beta_i}_*\left(\integ{\cD}{\cS}\right)_\pi$$
for $i=1,2$.
\end{claim}

Claim~\ref{indepclaim} immediately implies the sought independence on
the choices. Its proof will occupy us for the next several
subsections. Here we simply remark that the fine print in
\cite{MR2003c:14016} (specifically part (6) of Theorem 0.3.1) yields
the first part of the claim, as it guarantees that the inverse images
of the distinguished normal crossing divisors in $V_{\pi_i}$ are
normal crossing divisors; it is easily checked that all multiplicities
remain $>-1$ throughout the resolution.  What remains to be proved is
the claimed compatibility between the manifestations in
$V_\pi$ and $V_{\pi_i}$; by the recalled structure of the $\beta_i$'s,
it suffices to prove this in the particular case in which $\beta_1$,
$\beta_2$ are blow-ups along nonsingular centers meeting the relevant
divisors in $V_{\pi_1}$, $V_{\pi_2}$ with normal crossings.

\subsection{}\label{upshot}
The upshot of the preceding considerations is that the independence
of Definition~\ref{closedS} on the choices follows from the 
minimalist case of a particularly favorable blow-up, that is, the 
precise statement given below. Notation:
\begin{itemize}
\item $V$ is a nonsingular irreducible variety;
\item $E=\sum m_j E_j$ is a normal crossing divisor with nonsingular 
components $E_j$, $j\in J$, in $V$;
\item $\alpha: W \to V$ is the blow-up of $V$ along a nonsingular 
subvariety 
$B$ of codimension $d$, meeting $E$ with normal crossings;
\item $F_0$ is the exceptional divisor of the blow-up, and $F_j$ is the 
proper
transform of $E_j$, $j\in J$; let $J'=J\cup \{0\}$;
\item $m_0=(d-1)+\sum_{E_j\supset B} m_j$;
\end{itemize}

\begin{remark}
\begin{itemize}
\item Note that $F=\sum_{j\in J'} m_j F_j$ is a divisor with normal
crossings and nonsingular components (since $B$ meets $E$ with normal
crossings). Also, note that $F-\alpha^{-1}E=(d-1)F_0=K_\alpha$. This
is engineered to match the r\^ole of the normal crossing divisor
vis-a-vis the given divisor $\cD$ of the modification system,
cf.~\S\ref{intdef}: if $\cD$ is represented by $D_V$ on $V$ and $D_W$
on $W$, $\pi: V \to X$ is proper and birational, and 
$$D_V+K_\pi=\sum_{j\in J} m_j E_j$$
as in \S\ref{intdef}, then
$$D_W+K_{\pi\circ \alpha}=\alpha^{-1}(D_V+K_\pi)+K_\alpha=
\alpha^{-1}(\sum_{j\in J} m_j E_j)+(d-1) F_0=\sum_{j\in J'} m_j F_j$$
as needed in order for $W$ to again satisfy the assumptions given in 
\S\ref{intdef}.
\item The hypothesis that $B$ meets $E$ with normal crossings
implies that at most $d$ components $E_j$ contain $B$. Hence
$m_0>-1$ if all $m_j>-1$ for $E_j\supset B$, guaranteeing that the
assumption on multiplicities specified in \S\ref{intdef} is satisfied
in $W$ if and only if it is satisfied in $V$.
\end{itemize}
\end{remark}

Lastly, we must deal with $\cS$:
\begin{itemize}
\item $\Jbb$ is either the whole family of subsets of $J$, or the family 
of subsets of
$J$ having nonempty intersection with a fixed $L\subset J$;
\item $\Jbb'$ is the whole family of subsets of $J'$ in the first case; in 
the second 
case, $\Jbb'$ depends on whether any of the divisors $E_\ell$ for $\ell\in 
L$ 
contains the center $B$ of blow-up:
\item if none of the $E_\ell$ contains $B$, then $\Jbb'=\Jbb$;
\item if some of the $E_\ell$ contain $B$, then
$\Jbb'=\Jbb \cup \{I\,|\, I\subset J', 0\in I\}$.
\end{itemize}
This messy recipe encodes a rather simple situation. On $V$, $\cS$ is 
represented by either $V$ itself or by a union $\sum_{\ell\in L} E_\ell$, 
as required in \S\ref{intdef}. In the first case, $\cS$ is represented by 
$W$ on $W$; in the second case, it is represented by either the union 
$\sum_{\ell\in L} F_\ell$, if no $E_\ell$ contains $B$, or by 
$\sum_{\ell\in L\cup\{0\}} F_\ell$ if some $E_\ell$ do contain $B$. 
The prescription follows the fate of the distinguished families $\Jbb$, 
$\Jbb'$ through this predicament.

We are finally ready to state the main result, which will complete the
proof of Claim~\ref{indepclaim}, and hence of the independence of
Definition~\ref{closedS} on the choices.

\begin{theorem}\label{key}
With notation as above,
$$\alpha_*\left(
c(TW(-\log F))\cap \sum_{I\in \Jbb'}\prod_{i\in I}\frac{F_i}{1+m_i}
\right)
=c(TV(-\log E))\cap \sum_{I\in \Jbb}\prod_{i\in I}\frac{E_i}{1+m_i}
\quad.$$
\end{theorem}

\subsection{}
We prove Theorem~\ref{key} in \S\ref{Jbbisall} and \S\ref{Jbbisnotall}. 
The argument essentially amounts to careful bookkeeping, but is not 
completely straightforward. We collect a few necessary preliminaries in 
this subsection.

Let $V$ be a nonsingular variety, $B$ a nonsingular subvariety of 
codimension $d$, $\alpha: W\to V$ the blow-up of $V$ along $B$, and let 
$F$ be the exceptional divisor.
\begin{lemma}\label{necfacts}
The following hold in $A_*V$.
\begin{enumerate}
\item\label{log} Let $E$ be a divisor with normal crossings and 
nonsingular 
components $E_j$, $j\in J$, in $V$. Then
$$c(TV(-\log E))=\frac{c(TV)}{\prod_{j\in J} (1+E_j)}\quad.$$
\item\label{W} $\alpha_*(c(TW)\cap [W])=c(TV)\cap [V]+(d-1)\cdot c(TB)\cap 
[B]$.
\item\label{F} $\alpha_*(c(TF)\cap [F])=d\cdot c(TB)\cap [B]$.
\item\label{WF} $\displaystyle{\alpha_*\left(\frac{c(TW)}{(1+F)}\cap[W]
\right)=c(TV)\cap [V]-c(TB)\cap [B].}$
\item\label{factors} Let $E_j$, $j\in J$, be nonsingular hypersurfaces of 
$V$ meeting with normal crossings, and
let $F_j$ be the proper transform of $E_j$ in $W$. Assume at least one of 
the $E_j$ contains $B$. Then
$$\alpha_*\left(\frac{c(TW)}{(1+F)\prod_{j\in J}(1+F_j)}\cap [W]\right)=
\frac{c(TV)}{\prod_{j\in J}(1+E_j)}\cap [V]\quad.$$
\end{enumerate}
\end{lemma}

\begin{proof}
(\ref{log}) The equivalent statement $c(\Omega^1_V(\log E))$
$=c(\Omega^1_V)/\prod_{j\in J}(1-E_j)$ follows by a residue computation, 
see \cite{MR98d:32038}, 3.1.

(\ref{W}) This follows from Theorem~15.4 in \cite{MR85k:14004}.

(\ref{F}) Let $\underline\alpha$ be the projection $F\to B$, so that the 
class 
$\alpha_*(c(TF)\cap [F])$ is the push-forward to $A_*V$ of the class 
$\underline\alpha_*(c(TF)\cap [F])\in A_*B$.
The exceptional divisor is identified with the projectivization of the
normal bundle $N_BV$; therefore its tangent bundle fits in the Euler
sequence
$$\xymatrix{
0 \ar[r] & \mathcal O \ar[r] & \underline\alpha^* N_BV\otimes \mathcal 
O(1) \ar[r] & 
TF \ar[r] & \underline\alpha^* TB \ar[r] & 0
}\quad.$$
Hence
$$c(TF)=c((\underline\alpha^*N_BV\otimes \mathcal O(1))/\mathcal O)\cdot 
\underline\alpha^* c(TB)\quad,$$
and by the projection formula
$$\underline\alpha_*(c(TF)\cap [F])=c(TB)\cap \underline\alpha_*( 
c((\underline\alpha^*N_BV\otimes \mathcal O(1))/\mathcal O)\cap [F])
\quad.$$
Since $\underline\alpha$ has relative dimension $(d-1)$ over $B$, and 
$(\underline\alpha^*N_BV\otimes \mathcal O(1))/\mathcal O$ has 
rank~$(d-1)$,
only the top Chern class of this bundle survives the push-forward through 
$\underline\alpha$. This class may be evaluated using \cite{MR85k:14004},
Example~3.2.2, yielding the stated result.

(\ref{WF}) This can also be easily deduced from Theorem~15.4 in 
\cite{MR85k:14004}; or write
$$\frac{c(TW)}{1+F}\cap[W]=c(TW)\cap [W]-c(TF)\cap [F]$$
and apply (\ref{W}) and (\ref{F}).

(\ref{factors}) If $E_j$ does not contain $B$, then $F_j$ is the pull-back
of $E_j$; by the projection formula all such terms can be factored
out of both sides of the identity, so we may assume without loss of
generality that $E_j$ contains $B$ for {\em all\/} $j\in J$, with 
$J\ne\emptyset$. Let $J=\{1,\dots,r\}$, with $r\ge 1$.

The formula is immediate if $B$ is a hypersurface, so we may assume $B$
has codimension $>1$. Then $F_r$ is the blow-up of $E_r$ along $B$; and 
the other hypersurfaces $E_j$, $j<r$, cut out a divisor with normal
crossings and nonsingular components along $E_r$, containing $B$.
Also note that $\frac 1{1+F_r}=1-\frac {F_r}{1+F_r}$, and 
$c(TW)\frac {F_r} {1+F_r}\cap [W]=c(TF_r)\cap [F_r]$. If $r\ge 2$,
\begin{multline*}
\frac{c(TW)}{(1+F)\prod_{1\le j\le r}(1+ F_j)}
\cap[W]\\
=\frac{c(TW)}{(1+F)\prod_{1\le j<r}(1+ F_j)}\cap[W]
-\frac{c(TF_r)}{(1+F)\prod_{1\le j<r}(1+F_j)}
\cap[F_r]
\end{multline*}
and the needed equality follows if it is known for smaller, nonempty $J$.
Thus, we are reduced to proving
$$\alpha_*\left(
\frac{c(TW)}{(1+F)(1+\Til E)}\cap[W]\right)=
\frac{c(TV)}{1+ E}\cap[V]$$
for any nonsingular hypersurface $E$ containing the center $B$ of the
blow-up, where $\Til E$ denotes the proper transform of $E$. For this, 
rewrite 
the left-hand-side as
\begin{multline*}
c(TW)\left(1-\frac{F}{1+F}-\frac{\Til E}{1+\Til E}+\frac {F\Til E}
{(1+F)(1+\Til E)}\right)\cap [W] \\
=c(TW)\cap [W]-c(TF)\cap [F]-c(T\Til E)\cap [\Til E]+ c(T(F\cap \Til 
E))\cap 
[F\cap \Til E]
\end{multline*}
and use (\ref{W}) and (\ref{F}) to compute the push-forward:
\begin{multline*}
(c(TV)\cap [V]+(d-1)\cdot c(TB)\cap [B])-d\cdot c(TB)\cap [B] -
(c(TE)\cap [E]\\
+(d-2)\cdot c(TB)\cap [B])+ (d-1)\cdot c(TB)\cap [B]
=c(TV)\cap [V] - c(TE)\cap [E]\quad,
\end{multline*}
agreeing with the right-hand-side:
$$c(TV)\cap [V] - c(TE)\cap [E]=c(TV)\cap \left(1-\frac E{1+E}\right)\cap 
[V]=\frac{c(TV)}{1+E} \cap [V]$$
and concluding the proof.
\end{proof}

\begin{remark}
In characteristic~0, parts (\ref{W}) and (\ref{F}) are (even more) 
immediate from the functoriality of Chern-Schwartz-MacPherson classes.
\end{remark}

\subsection{}\label{Jbbisall}
In this subsection we prove Theorem~\ref{key} under the hypothesis that 
$\Jbb$ is the whole family of subsets of $J$. In this case, and using 
(\ref{log}) in Lemma~\ref{necfacts}, the statement to prove is
\begin{multline*}
\alpha_*\left(\frac{c(TW)}{(1+F_0)\prod_{j\in J}(1+F_j)}
\left(1+\frac{F_0}{M+d}\right)\prod_{j\in 
J}\left(1+\frac{F_j}{1+m_j}\right)
\cap [W]\right)\\
=\frac{c(TV)}{\prod_{j\in J}(1+E_j)}
\prod_{j\in J}\left(1+\frac{E_j}{1+m_j}\right)\cap [V]
\end{multline*}
with $M=\sum_{E_j\supset B} m_j$.

First, observe that if $E_j$ does {\em not\/} contain $B$, then 
$F_j=\alpha^*E_j$; by the projection formula, all factors involving such 
components can be factored out. Thus we may assume that all $E_j$ contain 
$B$, without loss of generality.

Second, with this additional assumption we can prove a stronger statement,
not binding $M$: we claim that, with $M$ an indeterminate,
\begin{multline*}\tag{*}
\alpha_*\left(\frac{c(TW)}{(1+F_0)\prod_{j\in J}(1+F_j)}
\left(1+\frac{F_0}{M+d}\right)\prod_{j\in 
J}\left(1+\frac{F_j}{1+m_j}\right)
\cap [W]\right)\\
=\frac{c(TV)}{\prod_{j\in J}(1+E_j)}
\prod_{j\in J}\left(1+\frac{E_j}{1+m_j}\right)\cap [V]
-\frac{M-\sum_{j\in J} m_j}{(M+d)\prod_{j\in J}(1+m_j)} c(TB)\cap [B]
\end{multline*}
if all $E_j$ contain $B$.

If $d=1$, that is, $B$ is itself a hypersurface of $V$, then $\alpha$ is
an isomorphism and (*) is immediately verified (note that in this case
$J$ consists of at most one element; and if $E_1$ contains $B$
then $F_0=B$, $F_1=0$ in $W\cong V$; the statement follows from the 
identity $\frac 1{M+1}=\frac 1{1+m_1}-\frac{M-m_1}{(M+1)(1+m_1)}$).

Therefore, we may assume $d\ge 2$.
Formula (*) is then proven by induction on the size of $J$.
If $J=\emptyset$, the statement is
$$\alpha_*\left(\frac{c(TW)}{(1+F_0)}
\left(1+\frac{F_0}{M+d}\right)\cap [W]\right)\\
=c(TV)\cap [V]-\frac{M}{M+d} c(TB)\cap [B]\quad.$$
To prove this, rewrite the left-hand-side as
$$\alpha_*\left(c(TW)\cdot\left(1-\frac{M+d-1}{M+d}\cdot \frac{F_0}{1+F_0}
\right)\cap [W]\right)\quad;$$
distributing and using (\ref{W}) and (\ref{F}) from Lemma~\ref{necfacts}:
$$c(TV)\cap[V]+(d-1)\cdot c(TB)\cap [B]- 
\frac{(M+d-1)d}{M+d} \cdot c(TB)\cap [B]$$
gives the stated result.

If $J=\{1,\dots,r\}$ with $r\ge 1$, split off the $F_r$ term from the 
left-hand-side 
of (*):
$$\frac{c(TW)}{(1+F_0)\prod_{j=1}^r(1+F_j)}
\left(1+\frac{F_0}{M+d}\right)\prod_{j=1}^r\left(1+\frac{F_j}{1+m_j}\right)
\cap [W]$$
equals
\begin{multline*}
\frac{c(TW)}{(1+F_0)\prod_{j=1}^{r-1}(1+F_j)}
\left(1+\frac{F_0}{M+d}\right)\prod_{j=1}^{r-1}\left(1+\frac{F_j}{1+m_j}\right)
\cap [W]\\
-\frac{m_r}{1+m_r}
\frac{c(TW)}{(1+F_0)\prod_{j=1}^{r-1}(1+F_j)}\frac{F_r}{1+F_r}
\left(1+\frac{F_0}{M+d}\right)\prod_{j=1}^{r-1}\left(1+\frac{F_j}{1+m_j}\right)
\cap [W]
\end{multline*}
Now note that $c(TW)\frac{F_r}{1+F_r}\cap [W]=c(TF_r)\cap [F_r]$;
that $F_r$ is the blow-up of $E_r$ along $B$ (since $E_r\supset B$,
and $d\ge 2$);
and that the other components cut $E_r$, $F_r$ along a divisor with
normal crossings. In other words, the induction hypothesis may be
applied to both summands in this expression. Care must be taken for
the r\^ole of $M$ in the second summand: as the codimension of $B$ in 
$F_r$ is $d-1$, the denominator $M+d$ must be viewed as $(M+1)+(d-1)$.
Therefore, applying the induction hypothesis evaluates the push-forward
as{\small
\begin{multline*}
\frac{c(TV)}{\prod_{j=1}^{r-1}(1+E_j)}
\prod_{j=1}^{r-1}\left(1+\frac{E_j}{1+m_j}\right)\cap [V]
-\frac{M-\sum_{j=1}^{r-1} m_j}{(M+d)\prod_{j=1}^{r-1}(1+m_j)} c(TB)\cap 
[B]\\
-\frac{m_r}{1+m_r}\frac{c(TE_r)}{\prod_{j=1}^{r-1}(1+E_j)}
\prod_{j=1}^{r-1}\left(1+\frac{E_j}{1+m_j}\right)\cap [E_r]
+\frac{m_r}{1+m_r}\frac{(M+1)-\sum_{j=1}^{r-1} m_j}{(M+d)
\prod_{j=1}^{r-1}(1+m_j)} c(TB)\cap [B]
\end{multline*}}
and now (*) follows by reabsorbing the $E_r$ term and performing
trivial algebraic manipulations:
\begin{multline*}
\frac{c(TV)}{\prod_{j=1}^{r-1}(1+E_j)}\left(1-\frac{m_r}{1+m_r} 
\frac{E_r}{1+E_r}
\right)\prod_{j=1}^{r-1}\left(1+\frac{E_j}{1+m_j}\right)\cap [V]\\
-\left(\frac{M-\sum_{j=1}^{r-1} m_j}{(M+d)\prod_{j=1}^{r-1}(1+m_j)}
-\frac{m_r}{1+m_r}\frac{M+1-\sum_{j=1}^{r-1} m_j}{(M+d)
\prod_{j=1}^{r-1}(1+m_j)}\right) c(TB)\cap [B]\\
=\frac{c(TV)}{\prod_{j=1}^r(1+E_j)}
\prod_{j=1}^r\left(1+\frac{E_j}{1+m_j}\right)\cap [V]
-\frac{M-\sum_{j=1}^r m_j}{(M+d)\prod_{j=1}^r(1+m_j)} c(TB)\cap [B]
\end{multline*}
as needed. This proves (*).

Setting $M=\sum_{E_j\supset B} m_j=\sum_{j\in J} m_j$ in (*)
concludes the proof of Theorem~\ref{key} for $\Jbb=$ the whole family 
of subsets of $J$.\qed

\subsection{}\label{Jbbisnotall}
Now assume that $\Jbb$ consists of the subsets of $J$ having nonempty 
intersection with a given $L\subset J$. Let $J=\{1,\dots, r\}$, and $L=\{1,
\dots,s\}$. As the proof of Theorem~\ref{key} in this case uses essentially
the same techniques as those employed in \S\ref{Jbbisall}, we provide
fewer details.

The statement to prove depends on whether some of the $E_\ell$, $\ell\in 
L$, contain $B$ or not. Using Lemma~\ref{necfacts} (\ref{log}), the
claim can be rewritten as follows:

($\dagger$) {\em If none of the $E_\ell$ contains $B$ for $\ell\in L$,\/} 
then 
\begin{multline*}
\alpha_*\left(
\frac{c(TW)}{(1+F_0)\prod_{j\in J}(1+F_j)} \left(1+\frac{F_0}{M+d}\right)
\sum_{\ell=1}^s\frac{F_\ell}{1+m_\ell}\prod_{i>\ell}\left(1+
\frac{F_i}{1+m_i}\right)\cap [W]
\right) \\
=\frac{c(TV)}{\prod_{j\in J}(1+E_j)} 
\sum_{\ell=1}^s\frac{E_\ell}{1+m_\ell}\prod_{i>\ell}\left(1+
\frac{E_i}{1+m_i}\right) \cap [V]
\end{multline*}

($\dagger\dagger$) {\em If some of the $E_\ell$ contain $B$ for $\ell\in 
L$,\/} 
then
\begin{multline*}
\alpha_*\left(
\frac{c(TW)}{(1+F_0)\prod_{j\in J}(1+F_j)}\frac{F_0}{M+d}
\prod_{j\in J}\left(1+\frac{F_j}{1+m_j}\right)\cap [W]
\right. \\
+\left.\frac{c(TW)}{(1+F_0)\prod_{j\in J}(1+F_j)}
\sum_{\ell=1}^s\frac{F_\ell}{1+m_\ell}\prod_{i>\ell}\left(1+
\frac{F_i}{1+m_i}\right)\cap [W]
\right) \\
=\frac{c(TV)}{\prod_{j\in J}(1+E_j)} 
\sum_{\ell=1}^s\frac{E_\ell}{1+m_\ell}\prod_{i>\ell}\left(1+
\frac{E_i}{1+m_i}\right) \cap [V]
\end{multline*}

\begin{proof}[Proof of $\text{\rm ($\dagger$)}$]
By the projection formula we may factor out all terms corresponding to
components not containing $B$; in particular, we may assume $L=\emptyset$, 
and the needed formula becomes
\begin{multline*}
\alpha_*\left(
\frac{c(TW)}{(1+F_0)\prod_{j\in J}(1+F_j)} \left(1+\frac{F_0}{M+d}\right)
\prod_{j\in J}\left(1+\frac{F_j}{1+m_j}\right)\cap [W]
\right) \\
=\frac{c(TV)}{\prod_{j\in J}(1+E_j)} 
\prod_{j\in J}\left(1+\frac{E_j}{1+m_j}\right) \cap [V]
\end{multline*}
with all $E_j$ containing $B$. This is precisely the formula proved in 
\S\ref{Jbbisall}.
\end{proof}

\begin{proof}[Proof of $\text{\rm ($\dagger\dagger$)}$]
By the projection formula we may, once more, assume that all $E_j$ 
contain~$B$.

Consider first the terms in the $\sum$ with $\ell\ge 2$, on the 
left-hand-side:
$$\frac{c(TW)}{(1+F_0)\prod_{j\in J}(1+F_j)}
\frac{F_\ell}{1+m_\ell}\prod_{i>\ell}\left(1+
\frac{F_i}{1+m_i}\right)\cap [W]\quad,\quad 2\le \ell\le s\quad.$$
We claim that each of these terms pushes forward to the corresponding
term in the $\sum$ on the right-hand-side:
$$\frac{c(TV)}{\prod_{j\in J}(1+E_j)} 
\frac{E_\ell}{1+m_\ell}\prod_{i>\ell}\left(1+
\frac{E_i}{1+m_i}\right) \cap [V]$$
(note: this is not so for the $\ell=1$ term!, cf.~Claim~\ref{explcom}).
To verify this, we argue as we did in \S\ref{Jbbisall}. The formula is 
clear if
$B$ has codimension~1; if $B$ has larger codimension, then $F_\ell$ is the 
blow-up of $E_\ell$ along $B$ (since all $E_j$ contain~$B$), and an 
induction on the number of factors reduces the verification to proving that
$$\alpha_*\left(\frac{c(TW)}{(1+F_0)\prod_{j\in J}(1+F_j)}\cap[W]\right)
=\frac{c(TV)}{\prod_{j\in J}(1+E_j)}  \cap [V]$$
if $J\ne\emptyset$ (this is where the hypothesis $\ell\ge 2$ is used).
This is proved in Lemma~\ref{necfacts},~(\ref{factors}).

Clearing the terms with $\ell\ge 2$ from both sides of $(\dagger\dagger)$,
we are reduced to proving that
\begin{multline*}
\alpha_*\left(
\frac{c(TW)}{(1+F_0)\prod_{j\in J}(1+F_j)}\frac{F_0}{M+d}
\prod_{j\in J}\left(1+\frac{F_j}{1+m_j}\right)\cap [W]
\right. \\
+\left.\frac{c(TW)}{(1+F_0)\prod_{j\in J}(1+F_j)}
\frac{F_1}{1+m_1}\prod_{i>1}\left(1+\frac{F_i}{1+m_i}\right)\cap [W]
\right) \\
=\frac{c(TV)}{\prod_{j\in J}(1+E_j)} 
\frac{E_1}{1+m_1}\prod_{i>1}\left(1+
\frac{E_i}{1+m_i}\right) \cap [V]
\end{multline*}
and this is implied by the following explicit computation:
\begin{claim}\label{explcom}
If all $E_j$ contain $B$, $j=1,\dots,r$, then:
\begin{multline*}
\alpha_*\left(\frac{c(TW)}{(1+F_0)\prod_{j=1}^r(1+F_j)}
\frac{F_1}{1+m_1}\prod_{i>1}\left(1+\frac{F_i}{1+m_i}\right)\cap [W]
\right) \\
=\frac{c(TV)}{\prod_{j=1}^r(1+E_j)} 
\frac{E_1}{1+m_1}\prod_{i>1}\left(1+
\frac{E_i}{1+m_i}\right) \cap [V]
-\frac{1}{\prod_{j=1}^r(1+m_j)}\cdot c(TB)\cap [B]\quad,
\end{multline*}
\begin{multline*}
\alpha_*\left(
\frac{c(TW)}{(1+F_0)\prod_{j=1}^r(1+F_j)}\frac{F_0}{M+d}
\prod_{j=1}^r\left(1+\frac{F_j}{1+m_j}\right)\cap [W]
\right) \\
=\frac{1}{\prod_{j=1}^r(1+m_j)}\cdot c(TB)\cap [B]\quad.
\end{multline*}
\end{claim}

\begin{proof}[Proof of the Claim]
The first formula is clear if $d=1$ (note that in this case $r=1$ 
necessarily, since the $E_j$ meet with normal crossings).

If $d\ge 1$ then each $F_j$ is the blow-up of $E_j$ along $B$. If $r>1$, 
splitting off the last factor and using
$$-\frac{1}{(1+m_1)\cdots (1+m_{r-1})}+\frac {m_r}{1+m_r} 
\frac{1}{(1+m_1)\cdots (1+m_{r-1})}
=-\frac{1}{(1+m_1)\cdots (1+m_r)}$$
shows that the general case follows from the case $r=1$:
\begin{multline*}
\alpha_*\left(\frac{c(TW)}{(1+F_0)(1+F_1)}
\frac{F_1}{1+m_1}\cap [W] \right) \\
=\frac{c(TV)}{(1+E_1)} 
\frac{E_1}{1+m_1}\cap [V]
-\frac{1}{(1+m_1)}\cdot c(TB)\cap [B]\quad.
\end{multline*}
Now this is equivalent to
$$\alpha_*\left(\frac{c(TF_1)}{(1+F_0)}
\cap [F_1] \right) \\
=c(TE_1)\cap [V]-c(TB)\cap [B]\quad,$$
which follows from (\ref{WF}) in Lemma~\ref{necfacts} (as $F_1$ is the 
blow-up of $E_1$ along $B$).

For the second formula, note that, for distinct $i_1,\dots,i_m$,
$$\alpha_*\left(\frac{c(TW)\cdot F_0\cdot F_{i_1}\cdots F_{i_m}}
{(1+F_0)\cdot(1+F_{i_1})\cdots (1+F_{i_m})}\cap [W]\right)
=(d-m)\cdot c(TB)\cap [B]\quad:$$
indeed, the intersection of the corresponding divisors $E_{i_1},\dots,
E_{i_m}$ is nonsingular, of dimension $\dim V-m$ (since the divisor 
$\sum E_j$ has normal crossings), and $F_0\cdot F_{i_1}\cdots 
F_{i_m}$ is the class of the exceptional divisor of its 
blow-up along $B$; so this follows from Lemma~\ref{necfacts}, (\ref{F}).

Therefore, rewriting the left-hand-side of the stated formula as
$$\alpha_*\left(
\frac 1{M+d}\frac{c(TW)\cdot F_0}{(1+F_0)}
\prod_{j=1}^r\left(1-\frac{m_j}{1+m_j}\frac{F_j}{1+F_j}\right)
\cap [W]\right)$$
and expanding the $\prod$ gives
$$\frac{1}{M+d}\left(d-(d-1)\sum \frac{m_j}{1+m_j}+(d-2) 
\sum \frac{m_jm_k}{(1+m_j)(1+m_k)}-\cdots\right) 
c(TB)\cap [B]$$
We leave to the reader the pleasant task of proving that this expression
equals the right-hand-side of the stated formula.\end{proof}

This concludes the proof of ($\dagger\dagger$), and hence of 
Theorem~\ref{key}.
\end{proof}

\subsection{}
At this stage the integral 
$$\integ{\cD}{\cS}\in A_*\cC_X$$ 
is defined for all {\em closed\/} subsets $\cS$ of $X$ and all
$\cD$. The definition for {\em locally closed\/} subsets is now forced
upon us: set
$$\integ{\cD}{\cS}:=\integ{\cD}{\cS_1}-\integ{\cD}{\cS_2}$$
if $\cS$ is the complement of a closed subset $\cS_1$ in a closed subset
$\cS_2$.

Of course we have to show that this is independent of the choices of the
closed subsets $\cS_i$: that is, if $\cS_i$, $\cT_i$ are closed in $\cC_X$,
and
$$\cS_1-\cS_2=\cT_1-\cT_2$$
as subsets of $\cC_X$, we have to prove that
$$\integ{\cD}{\cS_1}-\integ{\cD}{\cS_2}
=\integ{\cD}{\cT_1}-\integ{\cD}{\cT_2}\quad.$$
By comparing both sides with the intersection with $\cS_1\cap \cT_1$,
we may assume that $\cS_1\subset \cT_1$, and hence that
$\cT_1=\cS_1\cup \cT_2$ and $\cS_2=\cS_1\cap \cT_2$.
The needed equality is then the one arising from
$$(\cS_1\cup \cT_2)-\cT_2=\cS_1-(\cS_1\cap \cT_2)\quad,$$
that is, the following form of `inclusion-exclusion':

\begin{lemma}\label{incexc}
If $\cS$, $\cT$ are closed in $\cC_X$, and $\cD$ is a divisor, then
$$\integ{\cD}{\cS\cup\cT}=\integ{\cD}{\cS}+\integ{\cD}{\cT}
-\integ{\cD}{\cS\cap\cT}\quad.$$
\end{lemma}

\begin{proof}
The formula is clear if $\cS$ or $\cT$ equal $\cC_X$, so we
may assume that both are proper closed subsets.

It is enough to verify the statement for manifestations over a $\pi$
in $\cC_X$ resolving $\cD$ and $\cS$, $\cT$, $\cS\cap \cT$, and 
hence $\cS\cup \cT$. We may assume all are combinations of
components of a normal crossing divisor $\sum_{j\in J} E_j$: that is,
the subsets have ideals locally generated by $\prod_{j\in L_S}e_j$,
$\prod_{j\in L_T}e_j$, $\prod_{j\in L_S\cap L_T}e_j$, $\prod_{j\in
L_S\cup L_T}e_j$ respectively, where $e_j$ denotes a local generator
for $E_j$. Note that the ideal of $\cS\cap\cT$ in $V_\pi$ must be the
sum of the ideals of $\cS$ and $\cT$, to wit
$$(\prod_{j\in L_S}e_j)+(\prod_{j\in L_T}e_j)=
(\prod_{j\in L_S\cap L_T}e_j)\quad;$$
this implies that if $s\in L_S-L_T$ and $t\in L_T-L_S$, then the
corresponding components $E_s$, $E_t$ have empty intersection.
 
Now denote by $\Jbb_{\cS}$, $\Jbb_{\cT}$, $\Jbb_{\cS\cap \cT}$, 
$\Jbb_{\cS\cup \cT}$ the corresponding families of subsets $I\subset
J$. It suffices to prove that the individual contributions of each
$I\subset J$ to the expressions defining the integrals (as in
Definition~\ref{mandef}) satisfy the relation in the statement.

The only case in which this is not trivially true is when
$I\in\Jbb_{\cS}$ and $I\in\Jbb_{\cT}$, and hence $I\in \Jbb_{\cS\cup
\cT}$, but $I\not\in \Jbb_{\cS\cap \cT}$. In this case $I$ contains
indices $s\in L_S$, $t\in L_T$, such that $s\not\in L_T$, $s\not\in
L_S$; as observed above, this implies that $E_s\cap E_t=\emptyset$.
But then $\prod_{i\in I} E_i=0$, and hence the contribution of $I$
to all the integrals is zero. So the statement is verified in this
case as well, concluding the proof.
\end{proof}

\subsection{}\label{constr}
Finally, constructible sets are finite disjoint unions of locally 
closed subsets. If $\cS=\amalg_{k\in K} \cS_k$, with $\cS_k$ locally closed
in $\cC_X$, and $\cD$ is a divisor in $\cC_X$, we define
$$\integ{\cD}{\cS} := \sum_{k\in K} \integ{\cD}{\cS_k}\quad.$$
If $\cS_k$ is represented by $S_k$ in $V_\rho$, this element of
$A_*V_\rho$ is in fact the image of a class in $A_*(\cup_k \overline
S_k)$, cf.~Remark~\ref{support}. 

The integral $\integ{\cD}{\cS}$ ought to be independent of the ambient
system~$\cC_X$ containing $\cS$, but we only know instances of this
principle. For example, if $X$ is nonsingular and $\cS$ is represented
by $(id,S)$, with $S$ a subvariety of $X$, then the manifestation of
$\integ{0}{\cS}$ in $X$ resides naturally in $A_*S$, and is
independent of $X$, at least in characteristic~0: indeed, we will see
(\S\ref{CSMintro}) that it equals the Chern-Schwartz-MacPherson class
of $S$. Our argument will however rely on the theory of these classes;
we feel that there should be a more straightforward justification,
in the style of the computations performed in this section.

%%%

\section{Techniques of (celestial) integration}\label{properties}

\subsection{}\label{addit}
{\em Inclusion-exclusion.\/} 
The notion
$$\integ{\cD}{\cS}$$
introduced in \S\ref{integral} is additive in $\cS$ by definition
(cf.~\S\ref{constr}), and along the way we have had to prove a simple
form of inclusion-exclusion, that is, Lemma~\ref{incexc}. This implies
full inclusion-exclusion, that is: if $\cS=\cS_1\cap\cdots\cap \cS_r$,
then
$$\integ{\cD}{\cS}=\sum_{s=1}^r(-1)^{s-1}\sum_{i_1<\cdots<i_s}
\integ{\cD}{\cS_{i_1}\cup\cdots\cup \cS_{i_s}}\quad.$$
This is useful for explicit computations. For example, assume $\cS$ is
represented by $(X,S)$, with $S$ a subvariety of $X$; then $S$ may be
written as the intersection of hypersurfaces $S_1,\dots,S_r$, and
inclusion-exclusion reduces the computation of the integral over $\cS$
to the integral over constructible sets represented by hypersurfaces
of $X$.

\subsection{}\label{cov}
{\em Change-of-variables formula.\/} Let $\rho: Y\to X$ be a proper
birational map, with relative canonical divisor $K_\rho$. Here we need
that for proper birational maps $\pi: V \to Y$, with $V$ nonsingular,
$K_{\rho\circ\pi}=\pi^* K_\rho + K_\pi$; for example $Y$ could be
nonsingular, cf.~(\ref{funrcd}) in~\S\ref{relcandiv}.

Recall that $\cC_X$ and $\cC_Y$ are then equivalent, and
that consequently corresponding notions of divisors and constructible
subsets coincide, cf.~Remarks~\ref{remdivis} and~\ref{remconssub}; and
$A_*\cC_X\cong A_*\cC_Y$ by Lemma~\ref{isoeq}. So the integral of a
divisor $\cD$ over a constructible subset $\cS$ is defined in both
modification systems, and lands in the same target. 

\begin{theorem}
For all $\cS$ and $\cD$, and denoting by $K_\rho$ the divisor in
$\cC_X$ represented by $(\rho,K_\rho)$,
$$\integ{\cD}{\cS}=\integ[Y]{\cD+K_\rho}{\cS}\quad.$$
\end{theorem}

\begin{proof}
It is enough to prove the equality of manifestations in a variety $V$
resolving $\cD$ and $\cS$:
$$\xymatrix{
V \ar[r]^{\pi_Y} \ar@/_1pc/[rr]_{\pi_X} & Y \ar[r]^\rho & X
}$$
and this is clear from Definition~\ref{mandef}, since
$$\cD+K_{\pi_X}=(\cD+K_{\rho})+K_{\pi_Y}$$
as divisors in the modification system, by~(\ref{funrcd})
in~\S\ref{relcandiv}.
\end{proof}

\subsection{}\label{chernsimp}
{\em Chern classes.\/} If $X$ is nonsingular and $D$ is a nonsingular
hypersurface of $X$, let $\cD$ be the corresponding divisor of
$\cC_X$, represented by $(id,D)$. Then
$$\left(\integ{\cD}{\cC_X}\right)_{id}=c(TX)\cap [X]-\frac12
\cdot c(TD)\cap [D]\quad.$$
Indeed, the identity $id$ already resolves $\cD$, $\cC_X$; applying
Definition~\ref{mandef} gives
$$\left(\integ{\cD}{\cC_X}\right)_{id}=
\frac{c(TX)}{(1+D)}\cdot \left(1+\frac{D}2\right)\cap[X]=
c(TX)\cdot \left(1-\frac12\cdot \frac{D}{1+D}\right)\cap[X]$$
with the given result. In particular,
$$\left(\integ{0}{\cC_X}\right)_{id}=c(TX)\cap [X]\quad:$$
that is, the identity manifestation of the integral of $\one(0)$
realizes a `Chern class measure' on the variety. This observation is
extended readily to {\em nonsingular\/} subvarieties of $X$:

\begin{prop}\label{chernns}
Assume that $X$ is nonsingular, and $\cS\subset \cC_X$ is represented
by $(id,S)$. Then  
$$\left(\integ{0}{\cS}\right)_{id}=c(TS)\cap [S]\quad.$$
\end{prop}

\begin{proof}
Let $\pi:\Til X \to X$ be the blow-up of $X$ along $S$, with
exceptional divisor $E$. Then $\pi$ resolves $0$, $\cS$, and
$K_\pi=(d-1)E$, with $d=\codim S$. Applying Definition~\ref{mandef}
and Lemma~\ref{necfacts}, (\ref{log}) gives
$$\left(\integ{0}{\cS}\right)_{\pi}=
\frac{c(T\Til X)}{(1+E)}\cap \frac{E}d=\frac{c(TE)\cap [E]}d\quad;$$
from which the statement follows, by Lemma~\ref{necfacts}, (\ref{F}).
\end{proof}

We extend this remark to {\em singular\/} subvarieties in \S\ref{CSM}.

It is worth noting that the class in Proposition~\ref{chernns} is only
one manifestation of the integral. Manifestations in varieties mapping
to $X$ amount to specific lifts of the total Chern class of $X$. For
example, if $B$ is a nonsingular subvariety and $\pi:Y\to X$ is the
blow-up of $X$ along $B$, with exceptional divisor $E$, then
$$\left(\integ{0}{\cC_X}\right)_{\pi}=c(TY)\cap [Y]-\frac{d-1}d\cdot
c(TE)\cap [E]$$
according to Definition~\ref{mandef} (in one of the forms listed in
\S\ref{altexp}). Manifestations in other varieties are obtained as
push-forwards from the manifestation in a blow-up
(cf.~\S\ref{arbdef}). Explicit examples can be found in
Example~\ref{diffman}.

%%%

\section{Relation with Chern-Schwartz-MacPherson classes}\label{CSM}

\subsection{}\label{CSMintro}
Let $X$ be a nonsingular variety. In \S\ref{chernsimp} we have
observed that the symbol $d\mathfrak c_X$ behaves as a `Chern class 
measure' in the identity manifestation of the integral defined in
\S\ref{integral}, with respect to {\em nonsingular\/} subvarieties
$S$. That is:
$$\left(\integ{0}{\cS}\right)_{id}=c(TS)\cap [S]$$
if $\cS$ is the constructible subset of $\cS$ represented by 
$(id,S)$. It is natural to ask what class
$$\left(\integ{0}{\cS}\right)_{id}$$
computes if $S$ is not required to be nonsingular.

\begin{theorem}\label{CSMconn}
In characteristic~0
$$\left(\integ{0}{\cS}\right)_{id}=\csm(S)\in A_*X\quad,$$
the Chern-Schwartz-MacPherson class of $S$.
\end{theorem}

The Chern-Schwartz-MacPherson class generalizes in a beautifully
functorial fashion the Chern class of the tangent bundle of a
nonsingular variety. The reader is addressed to \cite{MR50:13587} for
the original work of Robert MacPherson (over $\Cbb$), and to
\cite{MR91h:14010} for a discussion over arbitrary algebraically
closed fields of characteristic~0. An equivalent notion had been
defined earlier by Marie-H\'el\`ene Schwartz, cf.~\cite{MR83h:32011}
for results comparing the two definitions.

The characteristic~0 restriction in the statement of
Theorem~\ref{CSMconn} is due to the fact that a good theory of
Chern-Schwartz-MacPherson classes does not seem to be available in
other contexts. Also, we have used the factorization theorem for
birational maps in the definition of the integral, and this relies on
resolution of singularities. In fact, we will work over $\Cbb$ in this 
section for simplicity of exposition, although everything can be
extended without difficulty to arbitrary algebraically closed fields
of characteristic~0.

Theorem~\ref{CSMconn} will follow from a more general result, linking the
integral to {\em MacPherson's natural transformation.\/}

\subsection{}\label{remind}
A quick reminder of related notions is in order. 

A {\em constructible function\/} on a variety~$X$ is a linear
combination of characteristic functions of closed subvarieties:
$\sum_{Z\subset X}=m_Z \one_Z$, where $m_Z\in \Zbb$ and $\one_Z(p)=1$
if $p\in Z$, $0$~otherwise. Thus a subset of $X$ is constructible if
and only if its characteristic function is.

Constructible functions form a group $F(X)$. Taking Euler
characteristic of fibers makes the assignment $X\to F(X)$ a covariant
functor under proper maps. More explicitly, if $f:X_1 \to X_2$ is a
proper map, and $Z$ is a subvariety of $X_1$, then the function
$f_*(\one_Z)$ defined by $p\mapsto \chi(f^{-1}(p)\cap Z)$ is
constructible; extending by linearity defines a push-forward $f_*:
F(X_1) \to F(X_2)$. 

MacPherson proved that there exists a unique natural transformation
$c_*$ from the functor $F$ to a homology theory---in this paper we use
Chow group with $\Qbb$-coefficients, denoted $A_*$, which also grants
us the luxury of using constructible functions with
$\Qbb$-coefficients---such that if $S$ is nonsingular, then
$c_*(\one_S)=c(TS)\cap [S]\in A_*S$. For arbitrary constructible
$S\subset X$, $c_*(\one_S)\in A_*X$ is the class denoted $\csm(S)$ in
\S\ref{CSMintro}.

An immediate application of the functoriality of the notion shows
that, if the ambient $X$ is complete, then the degree of $\csm(S)$
agrees with the topological Euler characteristic $\chi(S)$.

\subsection{}\label{IDSfun}
Let $X$ be a nonsingular variety over an algebraically closed field of
characteristic~0, and let $\cD$, $\cS$ resp.~be a divisor and a
constructible set in $\cC_X$. If $S$ is represented by $(\pi,S)$, and
$p\in X$, denote by $\cS_p$ the constructible set represented by
$(\pi, S\cap \pi^{-1}(p))$. 

\begin{defin}\label{IDSdef}
We define a function $I_X(\cD,\cS): X\to \Qbb$ by
$$I_X(\cD,\cS)(p):= \text{the degree of } \left(\integ{\cD}{\cS_p}
\right)_{id}\quad.$$
\end{defin}
\noindent (Recall that the manifestation on the right-hand-side may be
viewed as a class in $A_*p=\Qbb$, cf.~\S\ref{constr}, so its degree is
well-defined even if $X$ is not complete.)

The main result of this section is that the identity manifestation of
the integral defined in \S\ref{integral} corresponds to $I_X(\cD,\cS)$
via MacPherson's transformation:

\begin{theorem}\label{IDS}
The function $I_X(\cD,\cS)$ is constructible, and
$$\left(\integ{\cD}{\cS}\right)_{id}=c_*(I_X(\cD,\cS))\quad.$$
\end{theorem}

This result shows that, at the level of the identity manifestation,
the function $I_X(\cD,\cS)$ contains at least as much information as
the integral of $\one(\cD)$ over $\cS$. We could in fact define a
constructible function {\em of the system~$\cC_X$,\/} that is, an
element of the {\em inverse limit\/} of the groups of constructible
functions through the system (with push-forward as defined in
\S\ref{remind}), by setting
$$\cI_X(\cD,\cS)_\pi:=I_{V_\pi}(\cD+K_\pi,\cS)\quad.$$
Indeed, the following `change-of-variables' formula holds: if
$\pi'=\pi\circ \alpha$, then 
$$\alpha_* I_{V_{\pi'}}(\cD+K_{\pi'},\cS)=I_{V_\pi}(\cD+K_\pi,\cS)$$
(exercise!) The naturality of $c_*$ yields a homomorphism from the
group of constructible functions of the system to $A_*\cC_X$, and
the image of the `celestial' constructible function $\cI_X(\cD,\cS)$
through this homomorphism is the celestial integral.

Shoji Yokura has studied the {\em direct\/} limit of the groups of
constructible functions of an inverse system of varieties, and
associated Chern classes, in \cite{shoji}.

\subsection{}\label{IDSlemmas}
The proof of Theorem~\ref{IDS}, given in \S\ref{IDSproof}, relies on
two lemmas. We will use the following notation: if $E$ is a divisor
with components $E_j$, $j\in J$, and $I\subset J$, then $E_I^\circ$
denotes the complement of $\cup_{i\not\in I} E_i$ in $\cap_{i\in I}
E_i$.

\begin{lemma}\label{CSMlemma}
Let $V$ be a nonsingular variety, and let $E$ a divisor with normal
crossings and nonsingular components $E_j$, $j\in J$. Then for all
$I\subset J$
$$\csm(E_I^\circ)=\left(c(TV(-\log E))\cdot \prod_{i\in I} E_i\right) 
\cap [V]\quad.$$
\end{lemma}

\begin{proof}
Denote by $E_I$ the intersection $\cap_{i\in I} E_i$. Then $E_I$ is
nonsingular since $E$ has normal crossings, and its normal bundle in
$V$ has Chern class $\prod_{i\in I}(1+E_i)$, hence (using
Lemma~\ref{necfacts}, (\ref{log}))
$$\left(c(TV(-\log E))\cdot \prod_{i\in I} E_i\right) \cap [V]
=\frac{c(TE_I)}{\prod_{i\not\in I}(1+E_i)}\cap [E_I]\quad,$$
and we have to show that this equals
$$\csm(E_I^\circ)=c_*(\one_{E_I^\circ})=\csm(E_I)-\csm(E_I\cap(
\cup_{i\not\in I} E_i))$$
Now observe that $E_I\cap(\cup_{i\not\in I} E_i)$ is a divisor with
normal crossings in $E_I$; the needed formula follows then immediately
from (*) in \S2.2 of \cite{MR2001i:14009} (top of p.~4002).
\end{proof}

\begin{lemma}\label{denloe}
Let $V$, $W$ be nonsingular varieties; $\alpha: W \to V$ be a proper
birational map; $E$ a divisor with normal crossings and nonsingular
components $E_j$, $j\in J$, in $V$; $F$ a divisor with normal
crossings and nonsingular components $F_k$, $k\in K$, in $W$; and
$m_j$, $n_k$ integers such that 
$$\sum_{k\in K} n_k F_k=K_\alpha + \sum_{j\in J} m_j E_j\quad.$$
Finally, let $S$ be a constructible subset of $V$. Then
$$\sum_{I\subset J}\frac{\chi(E_I^\circ\cap S)}{\prod_{i\in I} (1+m_i)}
=\sum_{I\subset K}\frac{\chi(F_I^\circ\cap \alpha^{-1}(S))}
{\prod_{i\in I} (1+n_i)}
\quad.$$
\end{lemma}

\begin{proof}
By additivity of Euler characteristics we may assume that $S$ is
closed. The factorization theorem of \cite{MR2003c:14016} reduces the
statement to the case of a blow-up along a nonsingular center meeting
$E$ with normal crossings, which is worked out for the universal Euler
characteristic in Proposition~2.5 in \cite{math.AG/0401167}. 
\end{proof}

\begin{remark}\label{denloecon}
With $S=$ a point, and $E=0$, Lemma~\ref{denloe} states that
$$1=\sum_{I\subset K}\frac{\chi(F_I^\circ\cap 
\alpha^{-1}(p))}{\prod_{i\in I} (1+n_i)}\quad;$$
with $S=V$ and $E=0$ again, the statement is that 
$$\chi(V)=\sum_{I\subset K}\frac{\chi(F_I^\circ)}
{\prod_{i\in I} (1+n_i)}\quad.$$
These formulas have been known for a long time---they were first
proved by methods of $p$-adic integration in \cite{MR93g:11118}
(Theorem~6.1), and Fran\c cois Loeser informs me that he and Jan Denef
knew them as early as 1987; and that while aware of the implication
for Chern-Schwartz-MacPherson classes, they did not mention it for
lack of applications at the time.
These formula were later recovered (again by Denef and Loeser) by
using motivic integration (see for example the
survey~\cite{MR1905328}, \S4.4.3).
\end{remark}

In fact, this is the `point of contact with motivic integration'
mentioned in the introduction: 

\begin{claim}\label{contact}
Let $X$ be complete. With notations as in \S\ref{intdef}, and
$E^\circ_I$ as above,
$$\deg \integ{\cD}{\cC_X} = \sum_{I\subset J}\frac{\chi(E_I^\circ)}
{\prod_{i\in I}(1+m_i)}\quad.$$
\end{claim}

This follows immediately from the definition and from
Lemma~\ref{CSMlemma}, since the degree of $\csm$ equals the Euler
characteristic. This formula allows us to relate invariants
introduced by using our integral with other invariants arising
naturally from considerations in motivic (and/or $p$-adic)
integration, such as the {\em stringy Euler number,\/} cf.~\S7 of
\cite{math.AG/0401374}.

\subsection{}\label{IDSproof}
{\em Proof of Theorem~\ref{IDS}.}

Let $\pi:V \to X$ be an object of $\cC_X$ resolving $\cD$ and
$\cS$, cf.~\S\ref{intdef}. Thus there is a divisor $E$ with normal
crossings and nonsingular components $E_j$, $j\in J$, in $V$, such that 
$$D+K_\pi=\sum m_j E_j$$
and $(\pi,D)$ represents $\cD$; and $\cS$ is represented by $(\pi,S)$,
where $S=V$ or $S=\cup_{\ell\in L} E_\ell$ for some $L\subset J$.

First, we are going to show that
\begin{equation}\tag{$\star$}
I_X(\cD,\cS)=\pi_*\left(\sum_{I\subset J}\frac{\one_{E_I^\circ \cap S}} 
{\prod_{i\in I}(1+m_i)}\right)\quad,
\end{equation}
with push-forward of constructible functions defined as in
\S\ref{remind}; in particular, this shows that $I_X(\cD,\cS)$ is
constructible.

In order to show ($\star$), evaluate the right-hand-side at a $p\in X$:
$$\pi_*\left(\sum_{I\subset J}\frac{\one_{E_I^\circ \cap S}} 
{\prod_{i\in I}(1+m_i)}\right)(p)
=\sum_{I\subset J}\frac{\chi((E_I^\circ \cap S)\cap \pi^{-1}(p))} 
{\prod_{i\in I}(1+m_i)}$$
By Lemma~\ref{denloe}, this may be evaluated after replacing $\pi$
with an object dominating it and resolving $\cD$ and $\cS_p$
(represented by $S\cap \pi^{-1}(p)$, cf.~\S\ref{IDSfun}); that is, we
may assume that $S\cap \pi^{-1}(p)$ is a collection of components of
$E$, indexed by $L_p\subset J$. Let $\Jbb_{S_p}$ be the family of
subsets of $J$ meeting $L_p$. Then we may rewrite the right-hand-side
of ($\star$) as
$$\pi_*\left(\sum_{I\subset J}\frac{\one_{E_I^\circ \cap S}} 
{\prod_{i\in I}(1+m_i)}\right)(p)
=\sum_{I\in\Jbb_{S_p}}\frac{\chi(E_I^\circ)} 
{\prod_{i\in I}(1+m_i)}\quad.$$
Now using Lemma~\ref{CSMlemma} and the fact that the degree of the 
Chern-Schwartz-MacPherson class agrees with the Euler characteristic,
this equals the degree of
$$\sum_{I\in\Jbb_{S_p}}\frac{\csm(E_I^\circ)} {\prod_{i\in I}(1+m_i)}
=c(TV(-\log E))\cdot \left(\sum_{I\in\Jbb_{S_p}}\prod_{i\in I}
\frac{E_i} {1+m_i}\right)\cap [V]
=\left(\integ{\cD}{\cS_p}\right)_\pi\,.$$
Finally, the degree is preserved after push-forward, so this equals
$I_X(\cD,\cS)(p)$, concluding the proof of $(\star)$.

Now apply $c_*$ to both side of $(\star)$, and use
Lemma~\ref{CSMlemma} again. As in \S\ref{intdef}, denote by $\Jbb_S$
the family of subsets of $J$ if $S=V$, and the subfamily of subsets
meeting $L$ otherwise. This gives
\begin{align*}
c_*(I_X(\cD,\cS)) &= c_* \pi_*\left(\sum_{I\subset J}
\frac{\one_{E_I^\circ\cap S}} {\prod_{i\in I}(1+m_i)}\right)
= \pi_* \left(\sum_{I\subset J}\frac{c_*(\one_{E_I^\circ\cap S})} 
{\prod_{i\in I}(1+m_i)}\right)\\
&= \pi_* \left(\sum_{I\in \Jbb_S}\frac{\csm(E_I^\circ)} 
{\prod_{i\in I}(1+m_i)}\right)
= \pi_* \left(c(TV(-\log E))\cdot \sum_{I\in \Jbb_S}
\prod_{i\in I}\frac{E_i} {1+m_i}\right)\\
&= \pi_* \left(\integ{\cD}{\cS}\right)_\pi
= \left(\integ{\cD}{\cS}\right)_{id}\quad,
\end{align*}
concluding the proof of Theorem~\ref{IDS}.\qed

\subsection{}\label{claimproof}
Theorem~\ref{IDS} implies Theorem~\ref{CSMconn}. To see this, assume
$S\subset X$ is constructible. Then $I_X(0,\cS)(p)=0$ if $p\not\in S$,
since then $\cS_p=\emptyset$; if $p\in S$ and $X$ is nonsingular then 
$$I_X(0,\cS)(p)=\text{degree of }\left(\integ{0}{\cS_p}\right)_{id}=1$$
because then $\cS_p$ is represented by $(id,p)$, so the integral
computes $c(Tp)\cap [p]$, see~\S\ref{chernsimp}. This shows that
$I_X(0,\cS)=\one_S$ if $X$ is nonsingular, and hence
$$\left(\integ{0}{\cS}\right)_{id}=c_*(\one_S)=\csm(S)$$
by Theorem~\ref{IDS}.

Theorem~\ref{CSMconn} implies in particular that the identity
manifestation of the integral of $\one(0)$ over a subvariety $S$ of a
nonsingular $X$ is independent of the ambient variety~$X$.
We do not know to what extent integrals are independent of the ambient
variety, cf.~\S\ref{constr}.

\subsection{}
Given the close connection between the integral defined in
\S\ref{integral} and MacPherson's natural transformation, we feel that
a more thorough study of modification systems ought to yield a novel
approach to the theory of Chern-Schwartz-MacPherson classes. The
left-hand-side of the formula in Theorem~\ref{CSMconn} could be taken as
the {\em definition\/} of the class, and a good change-of-variable
formula for arbitrary proper maps should amount to the naturality of
this notion. With this (hypothetical) set-up, a proof of resolution of
singularity in positive characteristic would imply an automatic
upgrade of the theory of Chern-Schwartz-MacPherson classes in that
context.

In any case, if $S$ is any subvariety of a nonsingular variety $X$,
Theorem~\ref{CSMconn} affords many new manifestations of the
Chern-Schwartz-MacPherson class of $S$: for example, if $\pi:V \to X$
is any proper birational map, then
$$\left(\integ{0}{\cS}\right)_\pi$$
is a distinguished lift of $\csm(S)$ in $A_*V$. These manifestations
surely inherit good functoriality properties from $\csm(S)$, and it
would be interesting to explore these properties.

%%%

\newcommand{\cK}{{\mathcal K}}

\section{Applications}\label{applications}

\subsection{}
The mere existence of an integral satisfying the properties in
\S\ref{properties} has some immediate applications. For example,
assume that $X$ and $Y$ are nonsingular complete birational Calabi-Yau
varieties. Let $V$ be any resolution of indeterminacies of a
birational map between $X$ and $Y$:
$$\xymatrix{
& V \ar[dl]_{\pi_X} \ar[dr]^{\pi_Y}\\
X \ar@{-->}[rr] & & Y
}$$
Then $K_{\pi_X}=K_{\pi_Y}$, hence by change of variables (\S\ref{cov})
$$\integ[X]{0}{\cS}=\integ[V]{K_{\pi_X}}{\cS}=
\integ[V]{K_{\pi_Y}}{\cS}=\integ[Y]{0}{\cS}$$ 
for all $\cS$. Applying to $\cS=$ the whole system, and using that the
integral of $\one(0)$ evaluates the total Chern class
(\S\ref{chernsimp}), shows that the Chern classes of $X$ and $Y$ agree
{\em as elements of\/} $A_*\cC_X\cong A_*\cC_Y$. 

The same argument shows that any two nonsingular complete birational
varieties in the same $K$-equivalence class (cf.~Remark~\ref{remdivis})
have the same total Chern class in their (equivalent) modification
systems. This clarifies the main result of \cite{math.AG/0401167}.

\subsection{}\label{otherK}
The advantage of the more through investigation developed here over
the work in \cite{math.AG/0401167} is that we can now move away from
the hypothesis that the varieties are in the same $K$-class; in fact,
this can be done in several ways. For example, assume we have a
resolution of indeterminacies of a birational map between two
varieties $X$ and~$Y$, as above:
$$\xymatrix{
& V \ar[dl]_{\pi_X} \ar[dr]^{\pi_Y} \\
X \ar@{-->}[rr] & & Y
}$$
with $\pi_X$ and $\pi_Y$ proper and birational. The varieties $X$ and
$Y$ need not be complete, or nonsingular. The relative differentials
of $\pi_X$, $\pi_Y$ determine divisors $\cK_{\pi_X}$, $\cK_{\pi_Y}$ of
the equivalent modification systems $\cC_X$, $\cC_Y$. Let $\cD_X$,
$\cD_Y$ be any divisors such that $\cK_{\pi_X}+\cD_X=\cK_{\pi_Y}+\cD_Y$.

\begin{theorem}\label{spell}
With these notations, and for all constructible subsets $\cS$:
$$\integ[X]{\cD_X}{\cS}=\integ[Y]{\cD_Y}{\cS}\quad.$$
\end{theorem}

The proof is again an immediate application of change-of-variables.

Numerical consequences may be extracted from this formula. Recalling
that divisors act on Chow groups of modification systems
(Remark~\ref{remdivis}):

\begin{corol}\label{elgen}
With notations as above:
$$(c_1(X)-\cD_Y)^i\cdot\integ[X]{\cD_X}{\cS}
=(c_1(Y)-\cD_X)^i\cdot\integ[Y]{\cD_Y}{\cS}$$
for all $i\ge 0$.
\end{corol}

\begin{proof}
Indeed: $c_1(X)-\cD_Y=c_1(V)+\cK_{\pi_X}-\cD_Y=c_1(V)+\cK_{\pi_Y}-\cD_X
=c_1(Y)-\cD_X$.
\end{proof}

Specializing to the identity manifestation and taking degrees, in the
particular case in which $X$ and and $Y$ are nonsingular complete
varieties {\em in the same $K$-equivalence class,\/} taking
$\cD_X=\cD_Y=0$, and $\cS=$ the whole modification system, this gives
the equality:
\begin{equation}\tag{*}
c_1(X)^i\cdot c_{n-i}(X)=c_1(Y)^i\cdot c_{n-i}(Y)
\end{equation}
for all $i\ge 0$, with $n=\dim X$. 

For $i=0$ this is the well-known equality of Euler characteristics of
varieties in the same $K$-class (see for example \cite{MR2000i:14059}); 
for $i=1$ it can be derived as a consequence of the equality of
Hodge numbers, which follows from the change of variable formula in
motivic integration, and Theorem~3 in (\cite{MR91g:32039}).
The equality for all $i$ should also be a very particular case of the
fact that complex elliptic genera are preserved through
$K$-equivalence, cf.~\cite{MR2003j:14015}; and \cite{MR1953295}, where
this is byproduct of the definition of elliptic genera of singular
varieties (answering a fundamental question raised by Burt Totaro,
\cite{MR2001g:58037}, p.~758).

Even for varieties in the same $K$-equivalence class,
Corollary~\ref{elgen} is substantially stronger than $(*)$. For
example, let $S_X, S_Y$ be subvarieties of $X$, $Y$ resp., such that
$\pi_X^{-1}(S_X)=\pi_Y^{-1}(S_Y)$. Then
$$\deg(c_1(X)^i\cdot \csm(S_X))=\deg(c_1(Y)^i\cdot \csm(S_Y))\quad:$$
again take $\cD_X=\cD_Y=0$; and note that $S_X$, $S_Y$ represent the
same constructible subset in the modification system, then apply
Theorem~\ref{CSMconn}.

More generally, judicious choices for $\cD_X$, $\cD_Y$ may express
interesting data on $Y$ in terms of data on $X$; the simplest example
is probably
$$\chi(Y)=\deg\integ{\cK_{\pi_Y}-\cK_{\pi_X}}{\cC_X}\quad.$$
It is often possible to express $\cD_X$, $\cD_Y$ in terms of divisors
arising from subschemes of $X$, $Y$, and Corollary~\ref{elgen} may be
used to derive the equality of certain combinations of Chern numbers
of $X$, $Y$, and of these subschemes; in fact, Theorem~\ref{spell}
should simply be viewed as a notationally convenient way to encode a
large number of such identities. In general, these identities tend to
appear rather complicated, a lesson also learned through the work of
Lev Borisov and Anatoly Libgober, and Chin-Lung Wang. A very simple
prototypical situation is presented in Example~\ref{chernnos1}.

\subsection{}\label{birinv}
Some of the information exploited in \S\ref{otherK} is also captured
by the following invariant. For a nonsingular $X$, consider the set of
classes
$$\Can(X):=\left\{\integ{\cK}{\cC_X}\right\}\subset A_*\cC_X$$
as $\cK$ ranges over the divisors of $\cC_X$ obtained by pulling back
the effective canonical divisors of $X$.

{\em This is a birational invariant\/} of complete nonsingular
varieties. Indeed, so is $\Gamma(X,\Omega^{\dim X}(X))$, and the
change-of-variable formula ensures that the integrals of corresponding
divisors coincide. Taking degrees of the classes in $\Can(X)$ one
obtains a subset of $\Zbb$ which is likewise a birational invariant, 
and may be amenable to calculation. For example:

\begin{prop}
Let $X$ be a nonsingular complete algebraic variety, and assume that
$X$ is birational to a Calabi-Yau manifold $Y$. Then
$$\deg \Can(X)=\{\chi(Y)\}\quad.$$
\end{prop} 

\begin{proof}
Indeed $\Can(X)=\Can(Y)$ must be the single class
$\integ[Y]{0}{\cC_Y}$, whose identity manifestation is the total Chern
class of $Y$ by \S\ref{chernsimp}.
\end{proof}

\subsection{Zeta function}\label{zetafun}
For $\cD$ a divisor of $\cC_X$, and $m$ a variable, we can consider
the formal expression
$$Z(\cD,m):=\integ{m\cD}{\cC_X}\in A_*\cC_X[m]$$
(see \S\ref{negative} for a discussion of issues arising in letting
the multiplicities be variables). This `celestial zeta function' is a
very interesting object, which deserves further study. For example: 

\begin{prop}\label{topzetafun}
Assume $X$ is complete, and $\cD$ is the divisor corresponding to the
zero-scheme of a section $f$ of a line bundle on $X$. Then the degree
of $Z(\cD,m)$ equals the topological zeta function of $f$.
\end{prop}

We are referring here to the topological zeta function of
\cite{MR93g:11118}, see also \S6 in \cite{math.AG/0401374}, and
we are abusing the terminology since classically the topological zeta
function is defined for $f:M=\Cbb^n \to \Cbb$; this case can be
recovered by compactifying $M$ to $X=\Pbb^n$, then taking the integral
over the constructible subset represented by $(id,M)$).

Proposition~\ref{topzetafun} is proved easily by arguing as in
\S\ref{IDSproof} to relate $Z(\cD,m)$ to a combination of Euler
characteristics of subsets in the relevant normal crossing divisor in
a resolution, thereby matching the expression in \S6.6 of
\cite{math.AG/0401374}. 

The connection with the topological zeta function hints that the poles
of $Z(\cD,m)$ carry interesting information; a version of the {\em
monodromy conjecture\/} (\cite{math.AG/0401374}, \S6.8) can be phrased
in terms of $Z(\cD,m)$. For an explicit computation of $Z$, see
Example~\ref{zetaex}.

\subsection{Stringy Chern classes}\label{singclass}
The identity manifestation of the class
$$\integ{0}{\cC_X}$$ 
generalizes the total Chern class of the tangent bundle to possibly
singular~$X$. We loosely refer to this class as the {\em stringy\/}
Chern class of $X$, for reasons explained below.

There actually are different interpretations of this formula,
depending on the notion used to define the relative canonical divisor,
cf.~\S\ref{relcandiv}, and they lead to different classes. One
important alternative to the possibility presented in
\S\ref{relcandiv}, applicable to $\Qbb$-Gorenstein varieties, is to
let $\omega_X$ be the {\em double-dual\/} of the sheaf $\Omega^n_X$
(where $n=\dim X$). This is a divisorial sheaf, corresponding to a
Weil divisor $\Til K_X$; concretely, $\Til K_X$ may be realized as the
closure in $X$ of a canonical divisor of the nonsingular part of
$X$. The $\Qbb$-Gorenstein property amounts to the requirement that a
positive integer multiple $r\Til K_X$ of $\Til K_X$ is Cartier. If
$\pi: V \to X$ is a proper birational map, we can formally set $\Til
K_{\pi}$ to be a (`fractional') divisor such that $r\Til
K_\pi=rK_V-\pi^*(r\Til K_X)$. This definition satisfies the properties
mentioned in \S\ref{relcandiv}, hence it leads to an alternative
notion of integration in the modification system of~$X$.

In practice, the procedure sketched here assigns well-defined
multiplicities~$\in \Qbb$ to the components of the exceptional locus
of $\pi$, giving the input needed for the definition of the
integral. For more technical and contextual information on the
construction of $\omega_X$, see for example \cite{MR89b:14016}.

The class $(\integ{0}{\cC_X})_{id}$ obtained from this notion has some
right to be called the {\em stringy Chern class\/} of $X$, following
current trends in the literature (see e.g.~\cite{math.AG/0401374},
\S7.7); if $X$ is complete then the degree of its zero-dimensional
component equals the {\em stringy Euler number\/} of $X$, by
Claim~\ref{contact}. If $X$ admits a crepant resolution, the stringy
Chern class of $X$ is simply the image in $X$ of the Chern class of
any such resolution.

One difficulty with this notion is that it allows for the possibility that
some of the multiplicities $m_i$ in Definition~\ref{mandef} may be
$\le -1$, in which case our integral is simply undefined; this may
occur if the  singularities of $X$ are not log-terminal. This annoying
restriction may be circumvented in certain cases
(cf.~\S\ref{negative}) but appears to be necessary for the time being.

Example~\ref{compare} illustrates a clear-cut case in which the
classes obtained from the two notions of relative canonical divisor
considered here differ.

The choice of a notion for relative canonical divisors determines
a constructible function
$$I_X(0,\cC_X)$$
as in Definition~\ref{IDSdef}. By Theorem~\ref{IDS}, the corresponding
stringy class is the image via MacPherson's natural transformation of
this constructible function; which should hence be called the {\em
stringy\/} constructible function on $X$. (The analogous notion for a
Kawamata pair $(X,\Delta_X)$ would be the $\omega$ flavor of
$I_X(-\Delta_X,\cC_X)$.)

From our perspective these functions are more fundamental than their
incarnation as stringy Euler numbers or Chern classes: the information
carried by a stringy constructible function amounts to a list of
strata of $X$, and coefficients associated to these strata, from which
the invariants can be reconstructed by taking corresponding linear
combinations of the invariants of the strata. An alternative viewpoint
would associate to the stringy function a corresponding {\em stringy
characteristic cycle\/} in the cotangent bundle of a nonsingular
ambient variety containing~$X$. It is natural to guess that stringy
characteristic cycles admit a natural, intrinsic description.

It would be interesting to provide alternative computations of the
stringy functions (or characteristic cycles), possibly in terms
similar to those describing other invariants such as the local Euler
obstruction (which corresponds to the Chern-Mather class under
MacPherson's transformation). It would also be interesting to compare
the stringy class(es) to other notions of Chern classes for singular
varieties, such as Fulton's or Fulton-Johnson's
(cf.~\cite{MR85k:14004}, Example~4.2.6).

%%%

\section{Examples}\label{examples}
We include here a few explicit examples of computations of the
integral introduced in this paper.

\subsection{}\label{DintSex}
In \S\ref{chernsimp} we have seen that if $S\subset X$ are nonsingular,
and $\cS$ is represented by $(S,id)$, then
$\left(\integ{0}{\cS}\right)_{id}$ computes the Chern class of $S$.

More generally, if $\cD$ is represented by $(id,D)$, with $D$ a
nonsingular hypersurface intersecting $S$ transversally, then
$$\left(\integ{\cD}{\cS}\right)_{id}=c(TS)\cap [S]-\frac 12\cdot
c(T(D\cap S))\cap [D\cap S]\quad.$$
Indeed, the blow-up of $X$ along $S$ resolves $\cD$, $\cS$, and 
the formula is obtained easily from Definition~\ref{mandef} and
Lemma~\ref{necfacts}. Note that this shows that
$$\left(\integ{\cD}{\cS}\right)_{id}=
\left(\integ[S]{\cD_S}{\cC_S}\right)_{id}\quad,$$
where $\cD_S$ denotes the subset of $\cC_S$ represented by $(id,D\cap
S)$; that is, the left-hand-side is independent of the ambient
nonsingular variety $X$ in this case (cf.~\S\ref{constr}).

\subsection{}\label{IDSex}
An alternative (and more powerful) way to view the same computation is
through the use of the function~$I_X(\cD,\cS)$ introduced in
\S\ref{IDSfun}.

In the same situation ($X$ nonsingular, $\cD$ represented by $(id,D)$
with $D\subset X$ a nonsingular hypersurface), let $p\in X$; for fun,
consider a multiple $m\cD$ of $D$. If $p\not\in D$, then
$I_X(m\cD,\cC_X)(p)=I_X(0,\cC_X)=1$; if $p\in D$, by definition we can
compute $I_X(m\cD,\cC_X)(p)$ as the degree of 
$$\frac{c(T\Til X)}{(1+E)(1+\Til D)}\left(\frac{[E]}{m+n} +
\frac{[E\cup \Til D]}{(1+m)(m+n)}\right)$$
where $\Til X$ is the blow-up of $X$ along $p$, $E$ is the exceptional
divisor, $\Til D$ is the proper transform of $D$, and $n=\dim X$. A
straightforward computation evaluates this as $\frac 1{1+m}$, and
hence
$$I_X(m\cD,\cC_X)=\one_X-\frac{m}{1+m} \one_D$$
as should be expected.

Now if $\cS$ is represented by $(id,S)$, then
$$I_X(m\cD,\cS)=\one_S\cdot I_X(m\cD,\cC_X)=\one_S-\frac{m}{1+m}
\one_{S\cap D}\quad,$$
and by Theorem~\ref{IDS}
$$\left(\integ{m\cD}{\cS}\right)_{id}=\csm(S)-\frac m{1+m}\cdot
\csm(D\cap S)\quad,$$
generalizing the formula of Example~\ref{DintSex} to a constructible
$S$ with arbitrary singularities and intersecting $D$ as it wishes.

\subsection{}\label{Dnotdiv}
As we observed in Remark~\ref{remdivis}, every subscheme of $X$
determines a divisor in~$\cC_X$, which may be integrated. For example,
let $X$ be nonsingular, and let $Z\subset X$ be a nonsingular
subvariety of codimension~$d$; and let $\cZ$ be the divisor
corresponding to~$Z$, and $m\cZ$ the $m$-multiple of this divisor. Then
$$\left(\integ{m\cZ}{\cC_X}\right)_{id}=c(TX)\cap [X] 
- \frac m{d+m}\cdot c(TZ)\cap [Z]\quad.$$
Indeed, $\cZ=(\pi,E)$, where $\pi:\Til X \to X$ is the blow-up along
$Z$, and $E$ is the exceptional divisor; as $K_\pi=(d-1)E$,
Definition~\ref{mandef} gives the $\pi$ manifestation as
$$\frac{c(T\Til X)}{1+E}\cdot \left([\Til X]+\frac{[E]}{m+d}\right)
\quad,$$
from which the stated formula is straightforward (use
Lemma~\ref{necfacts}).

Of course there are divisors of $\cC_X$ which do not correspond to
subschemes of $X$. For example, suppose $D\subset X$ is a nonsingular
hypersurface containing a nonsingular subvariety~$Z$ of
codimension~$d$, and let $\cD$, $\cZ$ be the divisors of $\cC_X$
corresponding to $D$, $Z$. Then $\cD-\cZ$ is not represented in $X$,
even as a subscheme (if $d>1$). It is however represented by a divisor
in the blow-up along $Z$ (in fact, as the proper transform of~$D$) and
one computes easily that
$$\left(\integ{\cD-\cZ}{\cC_X}\right)_{id}=c(TX)\cap [X] 
-\frac 12 \cdot c(TD)\cap[D]+ \frac 1{2d}\cdot c(TZ)\cap [Z]\quad.$$

\subsection{}\label{diffman}
The Chern class of $\Pbb^2$ manifests itself in $\Pbb^1\times\Pbb^1$ as 
$$[\Pbb^1\times\Pbb^1]+\frac 32 [L_1]+\frac 32 [L_2] + 3[\Pbb^0]\quad,$$
where $L_1$ and $L_2$ denote lines in the two rulings.

This is obtained by considering the projection to $\Pbb^2$ from a
point~$p$ of a nonsingular quadric~$Q\cong\Pbb^1\times\Pbb^1$ in
$\Pbb^3$: the blow-up of $Q$ at $p$ resolves $0$, $\cC_{\Pbb^2}$, and
the computation in this blow-up is straightforward.

Similarly, the Chern class of $\Pbb^1\times\Pbb^1$ manifests itself in
$\Pbb^2$ as
$$[\Pbb^2]+\frac 52 [\Pbb^1] + 4[\Pbb^0]\quad.$$

The denominators in these two expressions imply the (otherwise
evident, in this case) fact that the birational isomorphisms between
$\Pbb^2$ and $\Pbb^1\times\Pbb^1$ does not extend to a regular map in
either direction.

\subsection{}\label{chernnos1}
A simple situation illustrating Theorem~\ref{spell} consists of a
birational morphism resolved by a blow-up and a blow-down:
$$\includegraphics[scale=.45]{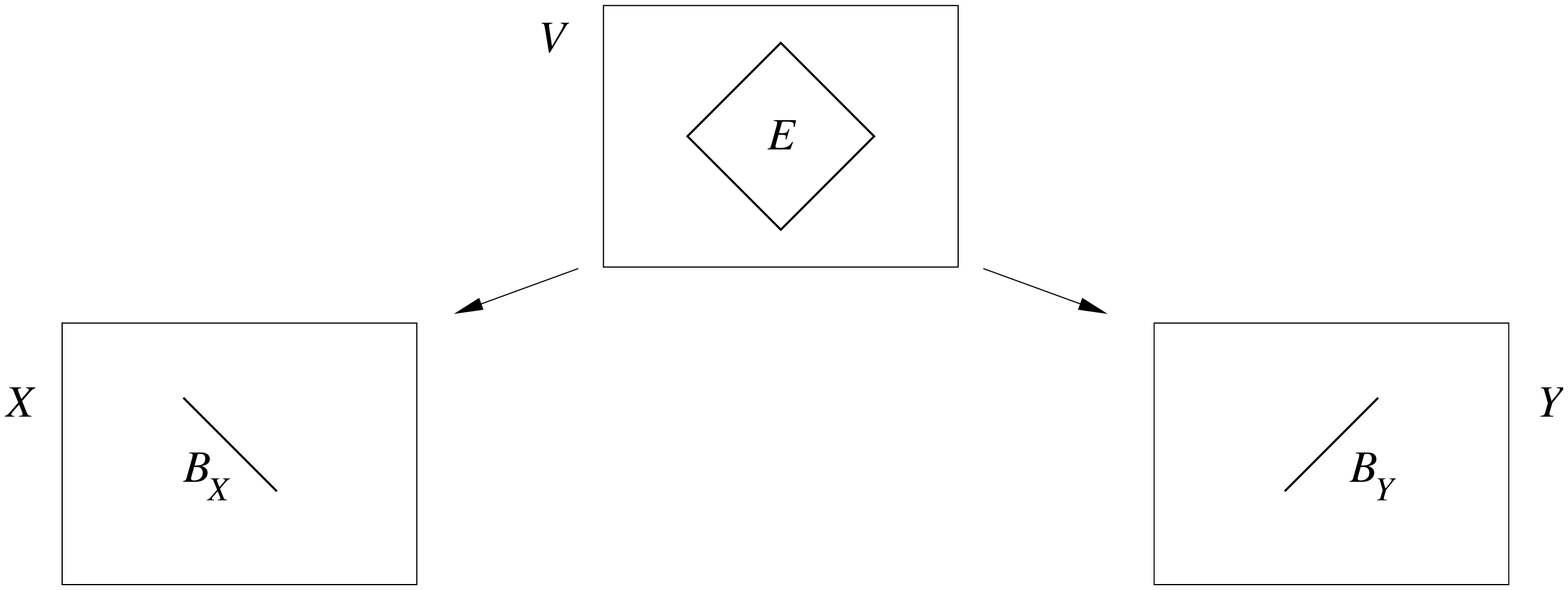}$$
Assume that $X$ and $Y$ are nonsingular and complete, $B_X\subset
X$, $B_Y\subset Y$ are nonsingular subvarieties, and $V=B\ell_{B_X}X =
B\ell_{B_Y}Y$, with $E$ the exceptional divisor for both blow-ups.
Let $d_X$, $d_Y$ be the codimension of $B_X$, $B_Y$ respectively.
By Theorem~\ref{spell}, and using the same notation for subvarieties
of $X$, $Y$ and for the corresponding divisors of the modification
systems,
$$\integ[X]{(d_Y-1)B_X+\cD}{\cC_X}=\integ[Y]{(d_X-1)B_Y+\cD}{\cC_Y}$$
for any divisor $\cD$. For example, representing $\cD$ by $(1-d_Y)E$
in the blow-up gives
$$\integ[X]{0}{\cC_X}=\integ[Y]{(d_X-d_Y)B_Y}{\cC_Y}\quad;$$
from which, evaluating the right-hand-side (using
Example~\ref{Dnotdiv}) and taking degrees:
$$\chi(X)=\chi(Y)+\frac{d_Y-d_X}{d_X}\chi(B_Y)\quad,$$
which is of course easy to check otherwise. With the same choice of
$\cD$, applying Corollary~\ref{elgen} with $i=1$ and taking degrees
gives a slightly more mysterious identity for other Chern numbers.
With $n=\dim X=\dim Y$, and denoting $c_{n-1}(X)$, etc.~for
$c_{n-1}(TX)\cap [X]$,~etc., one finds that
$$c_1(X)\cdot c_{n-1}(X)+(d_Y-d_X)\left(\frac{d_X+1}2 \chi(B_X) +
\frac 1{d_X} c_1(N_{B_X}X)\cdot c_{n-d-1}(B_X)\right)$$
must equal
$$c_1(Y)\cdot c_{n-1}(Y)+\frac{d_Y-d_X}{d_X}c_1(Y)\cdot
c_{r-1}(B_Y)\quad.$$
Making other choices for $\cD$, and varying~$i$, one easily gets a
large number of such identities.

If $d_X=d_Y$ the Chern numbers $c_1^i\cdot c_{n-i}$ for $X$ and $Y$
coincide ($X$ and $Y$ are in the same $K$-equivalence class, in this
case, cf.~\S\ref{otherK}).

If $d_Y=1$, then $V$ is isomorphic to $Y$, and the identities compare
Chern numbers of a variety $X$ and of its blow-up along a subvariety
$B_X$.

\subsection{Zeta function}\label{zetaex}
Let $X$ be a complete nonsingular surface, $D\subset X$ a nonsingular
curve, and $\cD$ represented by $(id,D)$. Then the identity
manifestation of $Z(\cD,m)$ is
$$\frac{c(TX)}{1+D}\left([X]+\frac {[D]}{1+m}\right)$$
and hence
$$\deg Z(\cD,m)=\chi(S)+\frac{m}{1+m}(K_S\cdot D+D^2)$$
where $K_S$, $\chi(S)$ are the canonical divisor and Euler
characteristic of $S$. By the adjunction formula, the `interesting
term' is a multiple of the Euler characteristic of~$D$.

Now assume that $D$ has a single singular point, consisting of an
ordinary cusp. The data $\cD$, $\cC_X$ is resolved by $\pi:\Til X\to
X$, obtained by a sequence of three blow-ups.
$$\includegraphics[scale=.6]{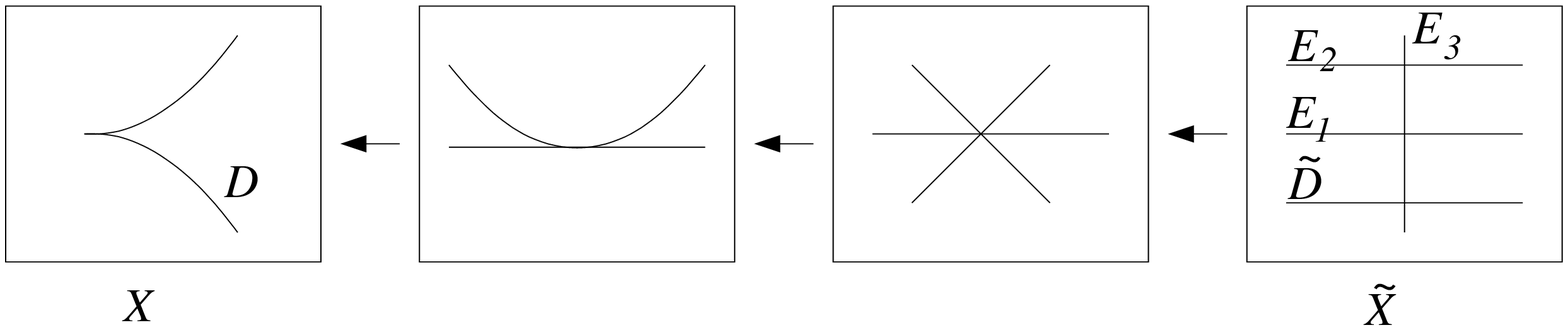}$$
The pull-back of $mD$ to $\Til X$ is the divisor
$$m\Til D+(1+2m)E_1+(2+3m)E_2+(4+6m)E_3\quad,$$
where $\Til D$ is the proper transform of $D$, and $E_1$, $E_2$, $E_3$
are (proper transforms of) the successive exceptional divisors.

Applying Definition~\ref{mandef} (in one of the equivalent forms in
\S\ref{altexp}), the $\pi$ manifestation of $Z(\cD,m)$ is
\begin{multline*}
c(T\Til X)\cap [\Til X]-\frac m{1+m}\cdot c(T\Til D)\cap [D]\\
-\frac{1+2m}{2+2m} c(TE_1)\cap [E_1]
-\frac{2+3m}{3+3m} c(TE_2)\cap [E_2]
-\frac{4+6m}{5+6m} c(TE_3)\cap [E_3]\\
+\left(\frac m{1+m}\cdot\frac{4+6m}{5+6m}
+\frac {1+2m}{2+2m}\cdot\frac{4+6m}{5+6m}
+\frac {2+3m}{3+3m}\cdot\frac{4+6m}{5+6m}\right)[p]
\end{multline*}
where $[p]$ is the class of a point. Taking degrees, we get
$$\deg Z(\cD,m)=\chi(S)+\frac m{1+m}(K_S\cdot
D+D^2)-\frac{12m}{5+6m}\quad.$$
Comparing with the nonsingular case, we can think of 
$$-\frac{12m}{5+6m}$$
as the contribution of the cusp to the zeta function. The fact that
this term has a pole at $-5/6$ is a trivial instance of the `monodromy
conjecture', see \S6.8 in \cite{math.AG/0401374}.

\subsection{Stringy classes}\label{compare}
We consider the poster example illustrating the distinction between
$\Omega^n_X$ and $\omega_X$ (\S\ref{relcandiv},~\S\ref{singclass}),
cf.~\cite{MR89b:14016}, \S1.8, and \cite{math.AG/0401374}, \S7.5.

Let $M$ be a nonsingular variety, and $X\subset M$ a subvariety with
the following property: there is a nonsingular subvariety $B$ of $X$
such that the proper transform $\Til X$ of $X$ in the blow-up $\Til M$
of $M$ along $B$ is nonsingular, and meets the exceptional divisor $E$
transversally. 
$$\xymatrix{
{\Til X} \ar[d]_\pi \ar[r]^j & {\Til M} \ar[d]^{\overline\pi} & E
\ar@{_(-}[l] \ar[d]\\
X \ar[r]^i & M & B \ar@{_(-}[l]
}$$

Let $n=\dim X$, and let $d$ be the codimension of $B$ in $X$; thus,
$d-1$ is the dimension of the (nonsingular) fibers of $\Til X\cap E$
over $X\cap B$.

\begin{claim}\label{normhyp}
($\Omega$ flavor.) With $K_\pi$ defined as in \S\ref{relcandiv},
$$K_\pi=(d-1)E\cdot \Til X\quad.$$
\end{claim}

To see this, note that $i^*\Omega^n_M \to \Omega^n_X$ is surjective,
hence the image of $\pi^*\Omega^n_X$ in $\Omega^n_{\Til X}$ is the
same as the image of $\pi^*i^*\Omega^n_M=j^*\overline\pi^*
\Omega^n_M$. At a point of $\Til X$ along $E$, by hypothesis we may
assume $\Til X$ has local coordinates $\Til x_1,\dots, \Til x_n$, part
of a system of coordinates of $\Til M$ (the other coordinates of $\Til
M$ being~$0$ along $\Til X$), mapping to $M$ according to
$x_1 = \Til x_1$, $x_2 = \Til x_2 \Til x_1, \dots, x_d = \Til x_d
\Til x_1$, $x_{d+1} = \Til x_{d+1},\dots, x_n = \Til x_n$. Here $\Til
x_1=0$ is the local equation of $E$ in $\Til M$. With these
coordinates, the only element of the evident basis of $\Omega^n_M$
surviving in $\Omega^n_{\Til X}$ is $dx_1\wedge\dots \wedge dx_n$,
which maps to $\Til x_1^{d-1}\, d\Til x_1\wedge\dots \wedge d\Til
x_n$. Thus the image of $\pi^*\Omega^n_X$ in $\Omega^n_{\Til X}$ is
$\Omega^n_{\Til X}\otimes \mathcal O(-(d-1)E)$, and the claim follows.

For the computation of the fancier $\Til K_{\pi}$, we will assume
further that $X$ is a hypersurface in $M$. Let $d, n$ be as above,
with $d\ge 2$, and let $k$ be the multiplicity of $X$ along $B$.

\begin{claim}\label{omegahyp}
($\omega$ flavor.) With $\Til K_\pi$ defined as in \S\ref{singclass},
$$\Til K_\pi=(d-k)E\cdot \Til X\quad.$$
\end{claim}

Indeed, as $\Til K_X$ agrees with the canonical divisor on the
nonsingular part of $X$, whose complement has codimension~$\ge 2$ in
$X$, the adjunction formula works as usual: $\Til K_X=(i^*K_M+X)\cdot
X$. Thus $\Til K_X$ is Cartier, and
$$\Til K_{\pi} = K_{\Til X}-j^*\pi^*(K_M+X) =K_{\Til X} - (K_{\Til 
M}-dE+\Til X+kE)\cdot \Til X=(d-k)E\cdot \Til X$$
(cf.~\cite{MR89b:14016}, p.~350).

Note that the multiplicity of this divisor is $\le -1$ for $k\ge d+1$;
this corresponds to the case in which the singularity is not
log-terminal. From our perspective, this limits computations of
integrals using this divisor to the case $k<(d+1)$
(cf.~\S\ref{negative}).

The two `stringy' Chern classes corresponding to the two flavors of
relative canonical divisors are respectively as follows:

\begin{prop}\label{stringyclasses}
\begin{align*}\tag{$\Omega$}
& \csm(X)+\frac 1d\cdot \left(\frac{(1-k)^{d+1}-1}k+1\right)
\cdot c(TB) \cap [B]\\
\tag{$\omega$}
& \csm(X)+\frac 1{d+1-k}\cdot \left(\frac{(1-k)^{d+1}-1}k+k\right)
\cdot c(TB) \cap [B] \qquad k< d+1
\end{align*}
\end{prop}

\begin{proof}
Since $\Til X$ resolves the data for the integral this is a
straightforward application of the definition, which we leave to the
reader. To obtain the given form (in terms of the
Chern-Schwartz-MacPherson class of~$X$), use Theorem~\ref{IDS}.
\end{proof}

The conventional stringy Euler number of $X$ is the degree of the
$\omega$ stringy Chern class, hence
$$\chi(X)+\frac 1{d+1-k}\cdot \left(\frac{(1-k)^{d+1}-1}k+k\right)
\cdot \chi(B)$$
by the second formula.

\begin{remark}
For this special class of hypersurfaces one can compute that
$$(1-(1-k)^{d+1})=\chi\quad,$$ 
the Euler characteristic of the Milnor fiber of $X$ at any of its
singularities, while $1-(1-k)^d = \text{Eu}$, the {\em local Euler
  obstruction\/} of \cite{MR50:13587}. It follows that
$$k=\frac{\text{Eu}-\chi}{\text{Eu}-1}\quad,$$
giving a more intrinsic flavor to the coefficients appearing in
Proposition~\ref{stringyclasses}. For example, the first class may be
rewritten
\begin{equation}\tag{$\Omega$}
\csm(X)+\frac 1d\cdot \frac{(\chi-1)\text{Eu}}{\chi-\text{Eu}}\cdot 
c(TB)\cap [B]\quad.
\end{equation}
\end{remark}

Taking degrees in Proposition~\ref{stringyclasses} gives explicit
formulas for the corresponding stringy Euler characteristics.

The second formula in Proposition~\ref{stringyclasses} shows that
there is no direct generalization of the $\omega$-stringy class to the
non-log-terminal case. On the other hand, the $\Omega$ flavor offers
an alternative which conveys essentially the same information (at
least for this simple-minded class of examples), and generalizes to
arbitrary singular varieties.

%%%

\section{Negative multiplicities}\label{negative}

\subsection{} 
According to the definition given in \S\ref{integral},
$\integ{\cD}{\cS}$ is simply {\em undefined\/} if the multiplicities
of the components of the relevant divisors in a resolving variety
happen to be $\le -1$. It would be desirable to weaken this
restriction, and we devote this final section to some musing on this
issue. This is not unimportant: for example, the invariant introduced
in \S\ref{birinv} is vacuous if the variety has no effective canonical
divisors; allowing integration over non-effective divisors would
enhance its scope. Also, allowing negative coefficients would extend
the range of the definition of stringy invariants to certain non-log
terminal singularities.

The key Definition~\ref{mandef}, that is, the manifestation of an
integral on a resolving variety, may be formally computed so long as
none of the multiplicities equals~$-1$; this would therefore appear to
be a more natural requirement than the stated one (that is, $m_j>-1$
for all $i$) for a resolving variety. The difficulty with adopting the
same definition with $m_j\neq -1$ is purely in the proof of
independence on the choice of resolving variety: while $m_j\neq -1$
may be satisfied for two chosen resolving varieties, it may fail at
some stage in the sequence of varieties connecting them through the
use of the factorization theorem. In other words, the first sentence
in Claim~\ref{indepclaim} does not hold as stated for this more
relaxed definition. This difficulty is also raised in a very similar
context by Willem Veys, \cite{math.AG/0401374}, \S8.1, Question~1.

One would like to show that if two varieties resolve the data $\cD$,
$\cS$ with all multiplicities $\ne -1$, then the resulting expressions
agree after push-forward. In particular, this would answers
affirmatively Veys' question. We note that allowing
multiplicities~$\le -1$ may cause some manifestations of the integral
to remain undefined; in such cases the integral would not exist as
an element of the Chow group $A_*\cC_X$, although it may still carry
useful information. The condition $m_j> -1$ adopted in
\S\ref{integral} avoids this problem in the simplest way.

\subsection{}
An example will highlight the difficulty, and will clarify how this
may be circumvented in certain cases.

\begin{example}\label{flop}
The quadric cone $X\subset\Pbb^4$ with equation $xy=zw$ is singular at
one point. This singularity may be resolved in several ways: blowing
up along the plane $x=z=0$ or along the plane $x=w=0$ produces two
`small' resolutions $X_-$, $X_+$ in which the singular point is
replaced by a $\Pbb^1$. Both these resolutions are dominated by the
blow-up $\Til X$ of $X$ at the vertex; in this blow-up the singularity
is replaced by a copy $E$ of $\Pbb^1\times\Pbb^1$. The birational
morphism from $X_-$ to $X_+$ is a classical example of {\em flop.\/}
$$\xymatrix{
& \Til X \ar[dl] \ar[dr] \ar[dd]_\pi\\
X_- \ar[dr]_{\pi_-} & & X_+ \ar[dl]^{\pi_+} \\
& X
}$$
The equation $x=0$ cuts out a Cartier divisor $D$ on $X$, consisting
of the union of the planes $x=z=0$, $x=w=0$ (these components are not
themselves Cartier divisors). We are interested in integrating the
divisor $-2\cD$ represented by $(id, -2D)$. 

The relative canonical divisor ($\omega$ flavor) is $0$ in both $X_-$
and $X_+$, since these are isomorphic to $X$ away from the distinguished
$\Pbb^1$'s. The divisor $D$ pulls back via $\pi_-$ to the union $D^-$
of two nonsingular divisors $D_1^-$, $D_2^-$ meeting transversally, so
we can formally apply Definition~\ref{mandef} and write
$$\left(\integ{-2\cD}{\cC_X}\right)_{\pi_-}=c(TX_-(\log D^-))\cap
\left( [X_-]+\frac{[D_1^-]}{1-2}+\frac{[D_2^-]}{1-2}+
\frac{[D_1^-\cap D_2^-]}{(1-2)^2}\right)$$
This can be evaluated easily, and pushes forward to
$$[X]+[D]$$
in $[X]$. An entirely analogous expression may be written in $X_+$,
with the same push-forward to $X$. This is as expected from
Claim~\ref{indepclaim}. 

However, the proof given in \S3 does not work in this case. The
blow-up/blow-down through $\Til X$ does give a factorization of the
flop. The inverse image $\pi^{-1}(D)$ consists of three components
meeting with normal crossings: two components $D_1$, $D_2$ dominating
the two components of $D$, and the exceptional divisor $E$. The
relative canonical divisor is $E$ (as shown in Claim~\ref{omegahyp}), so
$$-2\pi^{-1}(D)+\Til K_\pi=-2D_1-2D_2-E\quad.$$
Thus, $E$ appears with the forbidden multiplicity $-1$; the
corresponding formal application of Definition~\ref{mandef}:
\begin{multline*}
c(T\Til X(\log \pi^{-1} D))\cap
\left( [\Til X]+\frac{[D_1]}{1-2}+\frac{[D_2]}{1-2}+\frac{[E]}{1-1}
\right.\\
\left.+\frac{[D_1\cap D_2]}{(1-2)^2}+\frac{[D_1\cap E]}{(1-2)(1-1)}
+\frac{[D_2\cap E]}{(1-2)(1-1)}+\frac{[D_1\cap D_2\cap E]}{(1-2)^2
(1-1)}\right)
\end{multline*}
appears hopelessly nonsensical.
\end{example}

\subsection{}
We are now going to illustrate on this example how some information
may be extracted from such `meaningless' expressions in certain cases.
In an earlier version of this paper we in fact claimed that the same
approach used on this example can be applied in the general case; 
that appears to have been overly optimistic. I am grateful to Lev
Borisov and Wim Veys for pointing out problems with the original
argument. 

The idea is to view the multiplicity of $D$ as a variable $m$. The
basic formula guaranteeing independence on the resolution becomes an
equality of {\em rational functions\/} with coefficients in the Chow
group. The proof of Claim~\ref{indepclaim} goes through verbatim in
this context, and shows that the corresponding expression is well
defined as a rational function; therefore, so must be the expression
obtained by specializing $m$ to the needed multiplicity.

In Example~ref{flop}, the $\pi_-$ manifestation of the integral of
$m\cD$ is
$$\left(\integ{m\cD}{\cC_X}\right)_{\pi_-}=c(TX_-(\log D^-))\cap
\left( [X_-]+\frac{[D_1^-]}{1+m}+\frac{[D_2^-]}{1+m}+
\frac{[D_1^-\cap D_2^-]}{(1+m)^2}\right)$$
The $\pi$ manifestation makes sense as a rational function in $m$:
\begin{multline*}
c(T\Til X(\log \pi^{-1} D))\cap
\left( [\Til X]+\frac{[D_1]}{1+m}+\frac{[D_2]}{1+m}+\frac{[E]}{2+m}
\right.\\
\left.+\frac{[D_1\cap D_2]}{(1+m)^2}+\frac{[D_1\cap E]}{(1+m)(2+m)}
+\frac{[D_2\cap E]}{(1+m)(2+m)}+\frac{[D_1\cap D_2\cap E]}{(1+m)^2
(2+m)}\right)
\end{multline*}
and an explicit computation evaluates this as
\begin{multline*}
1+\frac{3+2m}{1+m}\,([D_1]+[D_2])+\frac{3+m}{2+m}\,[E]\\
+\frac{(2+m)(4+3m)}{(1+m)^2}\,[D_1\cap D_2]
-\frac{3+m}{1+m}\,[E]^2+\frac{(2+m)(3+m)}{(1+m)^2}\,[P]
\end{multline*}
where $[P]$ is the class of a point.

The many cancellations clearing several of the factors of $(2+m)$ at
denominator may appear surprising at first, but they are forced by the
strong constraints imposed on the situation by the geometry of the flop:
this class must push-forward {\em as a rational function\/} to the
$\pi_-$ and $\pi_+$ manifestations, which have no such factors; hence
only the terms that are collapsed by the push-forwards are allowed to
have poles at $m=-2$. In this example, the only such term is $E$ (as
$E^2$ and all other terms survive push-forward to $X_-$ and
$X_+$). Pushing forward to $X$ gives the identity manifestation:
$$\left(\integ{m\cD}{\cC_X}\right)_{id}=[X]+\frac{3+2m}{1+m}[D]
+\frac{(2+m)(4+3m)}{(1+m)^2}[L]+\frac{(3+m)(2+m)}{(1+m)^2}[p]\quad,$$
where $L$ is the class of a line. For $m=-2$, this agrees (of course)
with the class given in Example~\ref{flop}. For $m=0$ we obtain the
($\omega$ flavor of the) {\em stringy Chern class\/} of $X$, according
to the definition given in \S\ref{singclass} (and in agreement with
Proposition~\ref{stringyclasses}):
$$[X]+3[D]+8[L]+6[P]\quad.$$

As a last comment on this example, note that the situation is rather
different concerning the $\Omega$ flavor of the relative canonical
divisor: this is $2E$ for $\pi$ (by Claim~\ref{normhyp}), while it is
{\em not a divisor\/} for $\pi_-$ and $\pi_+$: for these maps,
$\Omega$ is tightly wrapped around the distinguished~$\Pbb^1$. In
particular, the $\Omega$ flavor of $\integ{m\cD}{\cC_X}$ cannot be
computed by applying Definition~\ref{mandef} to $X_-$ or $X_+$, as
these do not resolve the data. It may be computed by working in 
$\Til X$.

\subsection{}
A similar approach should work if, as in Example~\ref{flop}, the
key data comes from a divisor on $X$, for this guarantees that the
multiplicities on the two resolutions can be compatibly promoted to
variables, and Claim~\ref{indepclaim} shows that the corresponding
rational functions have the same push-forward in~$X$.
In fact, this attack to the question is not novel: see Remark~3.11 in
\cite{MR1953295} for a very similar approach; and the idea of
attaching variable terms to problematic multiplicities has been
championed by Veys with success in important cases
(\cite{MR2004h:14024}, \cite{MR2030094}).

We hope that the flexibility gained by considering the whole
modification system will allow us to assign compatible data in a more
general situation.

The conclusion to be drawn from the preceding considerations is that
the best setting in which to define the integrals considered in this
paper may in fact be a `decoration' of the Chow group of a
modification system $\cC_X$ by variables attached to all divisors of
$\cC_X$. The integral of a divisor should properly be considered as a
rational function with coefficients in the Chow group. The poles and
`residues' of such rational functions may carry valuable information,
as suggested by the connection with topological zeta functions
encountered in \S\ref{applications}.

%%%

%%%

\end{document}